\documentclass{llncs}
\usepackage{makeidx}  
\usepackage[pdftex]{graphicx}
\usepackage{amsmath}
\usepackage{amsfonts,amssymb}
\usepackage{clrscode3e}            
\usepackage{float}

\begin{document}
\frontmatter          
\pagestyle{headings}  
\addtocmark{Mixed-precision CG} 
\title{Parameter Optimization for Restarted Mixed-Precision Iterative Sparse Solver}
\titlerunning{Restarted mixed-precision iterative sparse solver}  
%
\author{Alexander Prolubnikov}
\authorrunning{Alexander Prolubnikov} 
%
\tocauthor{Alexander Prolubnikov}
\institute{Novosibirsk State University\\
\email{a.v.prolubnikov@mail.ru}
}
\maketitle              

\begin{abstract}

The problem of optimal precision switching for the conjugate gradient (CG) method applied to sparse linear systems is considered. A~sparse matrix is defined as an $n\!\times\!n$ matrix with $m\!=\!O(n)$ nonzero entries. The algorithm first computes an approximate solution in single precision with tolerance $\varepsilon_1$, then switches to double precision to refine the solution to the required stopping tolerance $\varepsilon_2$. Based on estimates of system matrix parameters --- computed in $O(m)$ time which does not exceed $1\%$ of the time needed to solve the system in double precision --- we determine the optimal value of $\varepsilon_1$ that minimizes total computation time. This value is obtained by classifying the input matrix using the $k$-nearest neighbors method on a small precomputed sample. Classification relies on a feature vector comprising: the matrix size $n$, the number of~nonzeros $m$, the pseudo-diameter of the matrix sparsity graph, and the~average rate of residual norm decay during the early CG iterations in~single precision. We show that, in addition to~the~matrix condition number, the diameter of the sparsity graph influences the growth of~rounding errors during iterative computations. The~proposed algorithm reduces the computational complexity of the CG --- expressed in equivalent double-precision iterations --- by more than $17\%$ on~average across the~considered matrix types in a sequential (non-parallel) setting. The~resulting speedup is at~most $1.5\%$ worse than that achieved with the optimal (oracle) choice of~$\varepsilon_1$.\\ 

While the impact of matrix structure on Krylov subspace method convergence is well understood, and the role of floating-point precision in~error accumulation has been extensively studied, the use of the sparsity graph diameter as a predictive feature for rounding error growth in mixed-precision CG appears to be novel. To the best of our knowledge, no prior work employs graph diameter to guide adaptive precision switching in~iterative linear solvers.

\end{abstract}

\section*{Introduction}

In mixed-precision algorithms, low-precision arithmetic is employed for efficiency, and higher-precision computations are used to refine the solution. For many problems in computational linear algebra, such algorithms can yield solutions with the same accuracy in less time than algorithms performing compu\-ta\-tions with fixed precision. Solving systems of linear algebraic equations, matrix factorizations, and finding eigenvalues and singular values are examples of such problems.

Work in this area has been ongoing since the early days of scientific and engineering computing, but most mixed-precision algorithms have been developed only in recent decades (see review papers \cite{Higham,Abdelfattah}). Traditionally, double-precision numbers were used for these computations, which today are 64-bit floating-point numbers. However, their use can be excessive and, in many cases, can be partially replaced by computations in single ($32$-bit) or even lower precision. On most modern hardware, single precision arithmetic is significantly faster than double-precision arithmetic. The reasons for this are that working with single-precision numbers requires half as much memory as double-precision numbers, reduces data transfer time, and lowers energy consumption. By switching precision in which operations are performed on a computer, we can use its resources more efficiently.

\smallskip

The precision with which real numbers are represented in a computer is characterized by the {\it machine epsilon}~$\epsilon$, defined as the distance between $1.0$ and the next larger floating-point number. The {\it unit roundoff} $u$, defined in IEEE 754 as $u\!=\!\epsilon/2$ for round-to-nearest arithmetic, is the maximum relative error incurred when rounding a real number to the nearest floating-point number: $$u=\max\limits_{x\in\mathbb{R},\ \mathrm{fl}(x)\neq 0}\frac{|x-\mbox{fl}(x)|}{|x|},$$ where $\mbox{fl(x)}$ denotes the floating-point representation of $x$. 

For the system of linear algebraic equations we are solving:
\begin{equation}\label{1}
Ax=b,
\end{equation}
the matrix $A\!\in\!\mathbb{R}^{n\times n}$ is symmetric, positive definite, and sparse. The matrix $A$ is called {\it sparse} if the number $m$ of nonzero elements in it is $O(n)$, and the constant in this estimate is such that the number of nonzero elements in it is no more than $10\%$ of the total number of its entries.

Let $x_i$ be the approximate solution of the system at the $i$-th iteration of~a~given method, and let $r_i\!=\!b-Ax_i$ be the {\it residual} of the approximate solution. We use the norm of the residual --- a {\it residual tolerance} or simply a {\it tolerance} --- as a convergence criterion. By the accuracy of the solution to the system of equations we mean the norm of the difference between the exact and approximate solutions. The accuracy of the solution depends on the conditioning of the matrix $A$.

The following approaches are mainly used to iterative methods of solving system of equations \cite{Higham}. The first, traditional approach to organizing mixed-precision computations is iterative refinement, in which the bulk of the computation is performed in low precision, while residuals or corrections are computed in higher precision. Let $u_1$ correspond to low precision, $u_2$ to high precision, $u$ to the precision with which the remaining operations of the computational algorithm are performed. The next computational scheme is used in the works of \cite{Turner}, \cite{Lindqwist} and many others. Let $i_{\max}$ specifies the maximum number of iterations in the scheme below. 

\begin{codebox}

\li Choose an initial $x_0$. 

\li  $i\!\leftarrow\!0$.

\li \Repeat

\li \Do
        Compute the residual $r_i\!=\!b\!-\!Ax_i$ in precision $u_2$.
				
\li     Solve the system $Ad_i\!=\!r_i$ in precision $u_1$.

\li     Compute $x_{i+1}\!=\!x_i+d_i$ in precision $u$.

\li     $i\!\leftarrow\!i+1$.
    \End
		
\li \Until $\bigl (i>i_{\max}\ \mbox{\bf or}\ \mbox{converged}\bigr )$				

\end{codebox}

Another mixed-precision approach is to use lower-precision computations to~find and/or apply the preconditioner. Thus, \cite{Giraud} proposes using a preconditioner computed in single precision for the conjugate gradient method, which is executed in double precision. Several variants of Jacobi preconditioning were considered, in which the diagonal blocks of the preconditioner are stored in half, single, or double precision for use with the conjugate gradient method. The precision is chosen according to a criterion based on an estimate of the condition number of the matrix block to which preconditioning is applied.

Mixed-precision computations in iterative methods can be employed for orthonormalizing the Krylov basis \cite{Gratton}, which is used in certain iterative solvers for linear systems. This idea is exploited in \cite{Aliaga}, where the authors propose storing the Krylov basis in reduced precision. Other approaches to using mixed-precision arithmetic in the solution of linear systems have also been proposed.

\smallskip

In the mixed-precision computation scheme under consideration, the first stage employs the conjugate gradient method (hereinafter referred to as CG) in single precision. In the second stage, once the residual norm of the approximate solution reaches a prescribed threshold, the computation switches to double precision, where the solution is refined until the required stopping tolerance is achieved. The residual tolerance at which this switch occurs is determined based on the matrix parameters. These parameters are computed in a time that is~small compared to the cost of running CG entirely in double precision.

This approach allows saving overall computation time by reducing the number of com\-pu\-ta\-tions in double precision. The savings achieved in this way, expressed in the number of double-precision iterations and taking into account the number of single-precision iterations performed, can reach $30\%$ or more.

\smallskip

The main problem that arises when applying the two-stage approach described above is the choice of the parameters of the matrix $A$ that influence the optimal choice of a tolerance of $\varepsilon_1$ achieved in the first stage. Even for matrices of the same type, the optimal values of $\varepsilon_1$ can differ significantly. To select $\varepsilon_1$ from a set of candidate values, we use the following parameters: the matrix size $n$, the number $m$ of nonzeros, the pseudo-diameter of the matrix sparsity graph (used as an estimate of its diameter), and the average rate of residual norm decay during the early iterations of CG in single precision. By the matrix sparsity graph (for brevity, we will henceforth refer to it simply as the {\it matrix graph}) we mean the undirected graph whose adjacency matrix is obtained by replacing all nonzero entries of $A$ with ones.

To determine the optimal value of $\varepsilon_1$ in terms of total computation time, we classify the input matrix using a feature vector composed of the parameters mentioned above. Each class is associated with one of the admissible values of~$\varepsilon_1$. Classification is performed using the $k$-nearest neighbors method. As will be shown below, the sample size required for its efficient operation is small, so~the time needed to classify an input matrix is significantly smaller than the cost of running CG entirely in double precision. Our experiments show that this approach is effective for a wide range of matrices, including random and banded ones.

\smallskip

The convergence behavior of Krylov subspace methods, such as CG method, is traditionally analyzed through spectral properties of the coefficient matrix --- most notably its condition number and eigenvalue distribution \cite{Saad2003,Greenbaum1997}. However, for~large sparse linear systems arising from the discretization of partial differential equations, the graph structure of the matrix also plays a crucial, albeit often implicit, role.

It is well understood that information propagates through the computational domain along the edges of the matrix graph during iterative solves \cite{Greenbaum1997,Simoncini2007}. The~diameter of this graph --- the maximum shortest-path distance between any two vertices --- governs the minimal number of iterations required for a disturbance (e.g., in the initial residual) to influence all components of the solution. For instance, in 2D and 3D grid-based problems, the graph diameter scales as $O(\sqrt{n})$ and $O(n^{1/3})$, respectively, which correlates with the growth of the condition number and, consequently, with CG convergence rates \cite{Saad2003}.

Although the graph diameter is rarely used as an explicit algorithmic parameter, related graph-theoretic concepts appear in the analysis of parallelism depth in incomplete factorizations \cite{Chow2015}. Furthermore, it is a well-established fact that the process of information propagation during iterative methods, and in particular during CG iterations, is governed by the matrix graph (see, e.g., \cite{Rude2026}). However, to the best of our knowledge, the graph diameter has not been previously employed as a predictive feature for adaptive precision switching in Krylov methods. In this work, we bridge this gap by demonstrating that the diameter---combined with other structural parameters---provides a robust indicator for determining the optimal switching point between single and double precision in CG, thereby minimizing total computational cost while maintaining solution accuracy.\\

\hfill

\newpage

\section{Two-Stage Algorithm for Solving System of Linear Equations with Mixed Precision}
We consider the following {\bf Algorithm I} for solving system (1); its scheme will be slightly modified below. The input consists of the matrix $A$ and right-hand side vector $b$, both given in double precision, the initial solution $x_0$ (taken to be zero), and the prescribed stopping tolerance $\varepsilon_2$ for the final refined approximate solution $\breve{x}$.

\begin{codebox}
\Procname{$\proc{\bf Algorithm I}$}
\li \mbox{\bf Input:}\ $A$, $b$, $x_0$, $\varepsilon_2$.\\

\li Determine the value of $\varepsilon_1$ for the matrix $A$.\\

\li Convert $A$ and $b$ to single precision: $\tilde{A}\!\leftarrow\!\texttt{float32}(A)$, $\tilde{b}\!\leftarrow\!\texttt{float32}(b)$.\\

\li Set the initial approximate solution to zero: $\tilde{x}_0\!\leftarrow\!0$.\\

\li Using CG, compute an approximate solution $\tilde{x}$ of the system $\tilde{A}x\!=\!\tilde{b}$\\ with tolerance $\varepsilon_1$.\\

\li Convert $\tilde{x}$ to double-precision vector $x_0$: $x_0\!\leftarrow\!\texttt{float64}(\tilde{x})$.\\

\li Using the initial approximate solution $x_0$, compute an approximate\\ solution $\breve{x}$ of the system of equations $Ax\!=\!b$ with stopping tolerance $\varepsilon_2$.\\

\li \mbox{\bf Output:}\ Refined approximate solution $\breve{x}$.

\end{codebox}

\noindent {\it The first stage} {\bf Algorithm I} is steps 1--4: determining the value of $\varepsilon_1$ for~the matrix $A$ and CG operation in single precision. {\it The second stage} is steps 5--6: CG~operation in double precision.\\

\section{Optimization Problem}
The problem of optimally choosing the parameter $\varepsilon_1$ is formulated as follows. The~parameter $\varepsilon_2$ specifies the required stopping tolerance for the iterative solution of the system. Given the matrix $A$, the goal is to choose $\varepsilon_1$ so as to minimize the total running time of {\bf Algorithm~I}. This is equivalent to minimizing an objective function of the form
\begin{equation}\label{3}
\hat{\varepsilon}_1(A)=\arg\min\limits_{\varepsilon_1} (t_1(A,\varepsilon_1)+t_2(A,\varepsilon_1,\varepsilon_2)),
\end{equation}
where $t_1$ is the execution time of the first stage of {\bf Algorithm I}, $t_2$ is the execution time of the second. The number of possible values of  $\varepsilon_1$ is small. In~the case we are considering, there are six of them: $$\varepsilon_1\!\in\!\{10^{-2},10^{-3},10^{-4},10^{-5},10^{-6},10^{-7}\}.$$

Let $N_1\!=\!N_1(A,\varepsilon_1)$ denote the number of iterations required to achieve residual tolerance $\varepsilon_1$ in the first stage, and $N_2\!=\!N_2(A,\varepsilon_1,\varepsilon_2)$ denotes the number of CG iterations performed in the second stage to reduce the residual norm below~$\varepsilon_2$. The computational complexity of {\bf Algorithm I} can be estimated as the~total number of iterations performed in both stages, which is expressed in~the~number of iterations performed in double precision: \begin{equation}\label{4}
I(A)=I(A,\varepsilon_1)=\omega N_1(A,\varepsilon_1)+N_2(A,\varepsilon_1,\varepsilon_2),
\end{equation}
\noindent where $\omega$ is the ratio of the average time $t_{\mathrm{sp}}$ required for a single-precision CG iteration in the first stage to the average time $t_{\mathrm{dp}}$ required for a double-precision CG iteration in the second stage: $\omega\!=\!t_{\mathrm{sp}}/t_{\mathrm{dp}}$.

Using criterion (3), we can formulate the optimization problem as follows:
\begin{equation}\label{4}
\hat{\varepsilon}_1(A)=\arg\min\limits_{\varepsilon_1} I(A,\varepsilon_1)=\arg\min\limits_{\varepsilon_1} \bigl (\omega N_1(A,\varepsilon_1)+N_2(A,\varepsilon_1,\varepsilon_2)\bigr ).
\end{equation}

We consider it more preferable to use (3) to estimate the running time of the algorithm when conducting experiments than using directly the total time $t_1+t_2$ of executing two stages on a computer. The processor time required by~a~computer to execute a program implementing an algorithm can be considered as a~random variable, generally speaking, with an unknown distribution. The value of this quantity is influenced by many factors that determine the operation of the computer. These factors include current processor load, data layout in memory, etc. By using (3), we are freed from the influence of these random factors when obtaining estimates of the computational complexity of the algorithms used and when analyzing the results of the experiments presented below. 

Among the factors not determined by the computer's technical implementation are the following. First, the vectors employed by the conjugate gradient (CG) method in single and double precision may exhibit different densities. Here, density refers to the proportion of nonzero entries in the vectors. This density influences both the computational cost per iteration in single and double precision and the resulting rounding errors. 

Secondly, the iteration runtime depends not only on the matrix itself and the floating-point precision used to represent its entries, but also on the right-hand side of the linear system. For different right-hand sides, both the execution time of a single CG iteration and the total number of iterations performed at each stage may vary. In our experiments, for each matrix in the sample, the right-hand side is generated randomly: we first draw a vector $x$ with components uniformly distributed over the interval $(1,2)$, and then compute $b\!=\!Ax$. 

Our estimates of the average execution time of CG iterations in single and double precision show that the ratio $\omega$ in (3)--(4) is typically much smaller than $1/2$, often reaching $1/5$ or lower in individual runs. This is consistent with results from the industrial problem-solving competition \cite{Huawei}, where a single-precision CG iteration was found to be $3$--$4$ times faster than in double precision.

In this regard, we use the value $\omega\!=\!1/3$. For the experiments we conducted, this choice almost always serves as an upper bound on the ratio $t_{\mathrm{sp}}/t_{\mathrm{dp}}$. Consequently, the efficiency estimates for $\omega\!=\!1/3$ presented below are pessimistic: in practice, the actual time savings are almost always greater. To illustrate this, we also provide estimates of the time savings obtained with a matrix-specific $\omega$, computed for each matrix in the sample by averaging the measured values of $t_{\mathrm{sp}}$ and $t_{\mathrm{dp}}$ over all CG runs performed on that matrix in single and double precision.\\

\section{Conjugate Gradient Method}
CG minimizes the quadratic function associated with the linear system (1): $$f(x)=\frac{1}{2} x^{\top}Ax-b^{\top}x.\eqno(5)$$ The minimum (5) is achieved by the exact solution of system (1).

\smallskip

CG is implemented by the algorithm given below. It receives as input the matrix $A$, the right-hand side $b$ of the system of equations, the initial guess $x_0$, and a required stopping tolerance of $\varepsilon$. The output of the algorithm is an approximate solution $\tilde{x}$ with required stopping residual tolerance.

\begin{codebox}
\Procname{$\proc{\bf CG}$}
\li \mbox{\bf Input:}\ $A$, $b$, $x_0$, $\varepsilon$.\\

\li $i\leftarrow 0$, $r_0 \leftarrow b-Ax_0$, $d_0\leftarrow r_0$.

\li\While $\|Ax_i-b\|>\varepsilon$:
   \Do 

\li $\alpha_i\leftarrow r_i^{\top}r_i/d_i^{\top}Ad_i$;

\li $x_{i+1}\leftarrow x_i+\alpha_id_i$;

\li $r_{i+1}\leftarrow r_i-\alpha_iAd_i$;

\li $\beta_{i+1}\leftarrow r_{i+1}^{\top}r_{i+1}/r_i^{\top}r_i$;

\li $d_{i+1}\leftarrow r_{i+1}+\beta_{i+1}d_i$;

\li $i\leftarrow i+1$;
    \End

\li $\tilde{x}\leftarrow x_i$.\\ 		

\li \mbox{\bf Output:}\ Approximate solution $\tilde{x}$.

\end{codebox}

The vectors $d_i$ are called \textit{search directions}, since at each iteration CG updates the approximate solution along this direction, which is a linear combination of~the~antigradient of~the~functional (5) --- i.e., the current residual vector --- and the previous search direction. 

In exact arithmetic, CG converges to the exact solution when the matrix $A$ is~symmetric positive definite. In finite-precision arithmetic, convergence is limited by rounding errors; however, if the floating-point precision is sufficiently high, CG can approach the solution to within machine precision before stagnating. The rate of convergence of CG is governed by the condition number $\kappa(A)$ of the system matrix: $$\kappa(A)=\|A\|\cdot\|A^{-1}\|=\frac{\lambda_{\max}}{\lambda_{\min}},$$ where $\lambda_{\min}$ and $\lambda_{\max}$ are the largest and smallest eigenvalues of $A$. The norm $\|\cdot\|_{A}$ of the error vector $e_i$ at the $i$-th iteration of CG ($\|e_i\|_{A}\!=\!(e_i^{\top}Ae_i)^{1/2}$) is~estimated as $$\|e_i\|_{A}\le 2\biggl (\frac{\sqrt{\kappa(A)}-1}{\sqrt{\kappa(A)}+1}\biggr )^i\|e_0\|_A.\eqno(6)$$ The number of iterations $k$ required to achieve $$\|e_k\|_A\!\le\!\varepsilon\|e_0\|_A$$ is bounded by $$k\le \left\lceil \frac{1}{2}\sqrt{\kappa(A)}\ln\biggl (\frac{2}{\varepsilon}\biggr )\right\rceil.\eqno(9)$$ 

The computational cost per iteration of CG is $O(m)$, where $m$ is the number of nonzeros in the sparse matrix $A$. The total number of iterations required to achieve a given accuracy typically scales with $\sqrt{\kappa(A)}$, leading to an overall complexity of $O\bigl (m\sqrt{\kappa(A)}\bigr )$ in the worst case. The memory complexity is $O(m)$.\\

\section{Parameters We Use to Determine $\varepsilon_1$}
Consider the parameters of the linear system matrix $A$ used to determine the value~of~$\varepsilon_1$.

Let $G\!=\!\langle V, E\rangle $ be the undirected, unweighted graph associated with the~matrix $A$, where $V\!=\!\{1,\ldots,n\}$ is the set of vertices corresponding to the rows and columns of $A$, and $(i,j)\!\in\!E$ if and only if $a_{ij}\!\neq\!0$. Since $A$ is symmetric positive definite, its diagonal entries are positive; thus, $G$ contains self-loops at every vertex. The graph $G$ is simple apart from these self-loops (i.e., it has no multiple edges), and for a~sparse matrix $A$, it is a {\it sparse graph}, meaning $m\!=\!|E|\!=\!O(n)$.

By a {\it quickly computable} parameter of the coefficient matrix of linear system (1), we mean a parameter that can be computed in a time not exceeding a~small fraction --- e.g., $1\%$ --- of the time required to solve the system using CG in~double precision. 

We consider the following parameters to determine the value of $\varepsilon_1$:
\begin{itemize}
\item[1)] the size $n$ of a square matrix $A$;
\item[2)] the number of its nonzero entries $m$ (i.e., the number of edges in the graph of the matrix $A$); 
\item[3)] pseudo-diameter $\tilde{\ell}(A)$ of the matrix graph --- estimate of the diameter of the matrix graph $A$;
\item[4)] the average rate $v$ of reduction in the CG residual norm during the early single-precision iterations.
\end{itemize}

Thus, to predict the optimal value $\hat{\varepsilon}_1$ given by (4), we use the following parameter vector: $$\chi(A)=(n,m,\tilde{\ell},v).$$

The matrix size $n$ and the number of nonzeros $m$ influence the number of CG iterations required at any floating-point precision, so we include them in the parameter vector $\chi(A)$ used to predict $\hat{\varepsilon}_1$. The value of $n$ determines the maximum number of CG iterations needed to solve the system in exact arithmetic. The value of $m$ characterizes the matrix density and governs the cost of the matrix-vector multiplication --- the dominant operation in CG. Given finite floating-point precision, both parameters affect the cost per iteration, the total iteration count, and the accumulation of rounding errors; consequently, they must be accounted for when predicting $\hat{\varepsilon}_1$. Since $n$ and $m$ are obtained while reading the matrix $A$ from file, they incur no additional computational overhead.

\smallskip

Consider the parameters $\tilde{\ell}$ and $v$, which we also use to predict $\hat{\varepsilon}_1$ for the matrix $A$, and compare the computational cost of computing these parameters and performing the prediction with that of solving the linear system using CG entirely in double precision.\\ 

\section{Growth of Rounding Errors During Iterations\\ and the Diameter of the Matrix Graph}

\subsection{Topological properties of a matrix}

Typically, in the analysis of CG convergence---estimating convergence rates, residuals, and errors at its iterations---one employs {\it spectral} properties of the matrix: the condition number (i.e., the extreme eigenvalues) and the overall distribution of the matrix spectrum. When spectral parameters are used, the matrix is treated as a linear operator. Under this approach, when analyzing CG convergence for arbitrary matrices, it suffices to consider diagonal matrices \cite{Carson2024}.

When working with {\it topological} parameters of a matrix, the coefficient matrix of the linear system is interpreted as the adjacency (or structure) matrix of a graph that encodes the connectivity pattern defined by the matrix---\textit{sparcity pattern} of the matrix. As will be shown below, such topological characteristics---namely, distances between vertices and the diameter of the matrix graph--- strongly influences the accumulation of rounding errors incurred during the execution of CG.

\smallskip

The diameter of the matrix graph is related to the spectral properties of the matrix. In particular, the number of distinct eigenvalues is at least the diameter of the graph \cite{VanDamHaemers1995}. Typically, matrices whose graphs have a large diameter are ill-conditioned, whereas, for example, in such practically important classes of sparse matrices as M-matrices, symmetric diagonally dominant matrices, and grid-based matrices, a small diameter is more often associated with good conditioning---although a small diameter by itself does not guarantee good conditioning.

Considering the diameter of the matrix graph as a parameter influencing the accumulation of rounding errors during CG iterations is justified for sparse matrices, whose graph diameters can differ significantly. In contrast, for dense matrices with $O(n^2)$ nonzeros, the sparcity pstterns corresponding to matrices of the same dimension $n$ are nearly identical; consequently, the differences in the number of potential rounding errors which are caused by topological properties we consider will also be negligible.

\subsection{Rounding error propagation during iterative process $b_i=Ab_{i-1}$}

It is well known that the process of information propagation during CG iterations is governed by the matrix graph (see, e.g., the work \cite{Rude2026}, where this process is vividly illustrated). By \emph{information propagation} during the execution of an iterative method for solving a linear system, we mean the spread of updates to the components of the approximate solution vector throughout the iterations. Updated components---computed with rounding errors---are subsequently used to compute other components. Consequently, rounding errors introduced in earlier iterations propagate together with the updated solution components during the current iteration.

Let $G\!=\!\langle V, E\rangle $ be the {\it matrix graph}---an undirected, unweighted graph---associated with the matrix $A$, where $V\!=\!\{1, \ldots, n\}$ is the set of vertices corresponding to the rows (and columns) of $A$. An edge $(i,j)\!\in\!E$ exists if and only if $a_{ij}\!\ne\!0$. The graph $G$ is simple (i.e., undirected and unweighted) but includes self-loops at every vertex, since the diagonal entries of a symmetric positive definite matrix $A$ are strictly positive. If $A$ is sparse, then its graph $G$ is also {\it sparse}, meaning that the number of edges satisfies $m\!=\!|E|\!=\! O(n)$.

A {\it path} in a graph is a finite sequence of distinct vertices such that each vertex, except the last one, is connected to the next by an edge. For  a {\it walk}, vertices and edges may repeat. The {\it length} of a path or a walk is the number of edges connecting the vertices in the sequence. Let $l(i,j)$ denote the distance between vertices $i, j\!\in\! V(G)$ in the graph $G$ associated with the matrix $A$, defined as the number of edges in the shortest path from $i$ to $j$. Since the matrix is symmetric, the graph is undirected, and thus $l(i,j)\!=\!l(j,i)$. The {\it diameter} $\ell$ is the greatest distance between any pair of vertices in the graph: $$\ell\!=\!\max_{i,j \in V} l(i,j).$$ 

Consider the following iterative process: $$b_1 = A b_0,\quad b_2 = A b_1,\quad \ldots,\quad b_k = A b_{k-1}.\eqno(8)$$ Let $\mathcal{N}_k\!=\!\{(i,j) \mid l(i,j)\!=\!k\}$ denote the set of unordered pairs of vertices $(i,j)$ in the graph $G$ that there are walks of length $k$ connected them (they may include loops). Since the $k$-th power of the matrix $A$ has a nonzero entry $(A^k)_{ij}$ if and only if there exists a walk of length $k$ between vertices $i$ and $j$, the first appearance of a connection between $i$ and $j$ occurs precisely at $k\!=\!l(i,j)$, it follows that the number of vertex pairs connected by walks of length $k$ is reflected in the sparsity pattern of $A^k$.

When computing the matrix-vector product $A\cdot b_i$, each traversal of an edge $(i,j)\!\in\!E(G)$ corresponds to multiplying the matrix entry $a_{ij}$ by the component $b_{ij}=(b_i)_j$ of the vector $b_i$ and adding the result to the current value of the $i$-th component being updated. Each such arithmetic operation introduces a potential source of rounding error.

Since $$b_{1j}\!=\!\sum_{i=1}^n a_{ji} b_{0i}, \eqno(9)$$ the computation of a single component $b_{1j}$ involves multiplications of the form $a_{ji} \cdot b_{0i}$ for all indices $i$ such that $a_{ji}\ne 0 $, i.e., for all pairs $(j,i)\!\in\!\mathcal{N}_1$. For each such multiplication, one floating-point addition is subsequently performed to accumulate the partial sum in (9). Hence, for every nonzero entry $a_{ji}$, two elementary floating-point operations are executed: one multiplication and one addition---each potentially introducing a rounding error.

Furthermore, $$b_{2j} = \sum_{i=1}^n a_{ji} b_{1i}, \eqno(10)$$ so if $a_{ji}\!\ne\!0$ and the component $\!b_{1i}\!$ was updated during the first iteration---that is, if its value differs from the original $b_{0i}$---then this change is propagated to the component $b_{2j}$.  

Consequently, the updates applied to the initial vector $b_0$ to obtain $b_1$ are further propagated to those vertices $t\!\in\!V(G)$ for which there exists an intermediate vertex $j$ such that $$(a_{ij}\ne 0) \land (a_{jt}\ne 0), \eqno(11)$$ i.e., to all vertex pairs $(i,t)$ connected by a walk of length two. In other words, these are precisely the pairs satisfying $(i,t)\!\in\!\mathcal{N}_2$.


\begin{figure}[H]
\centering
\includegraphics[width=25mm]{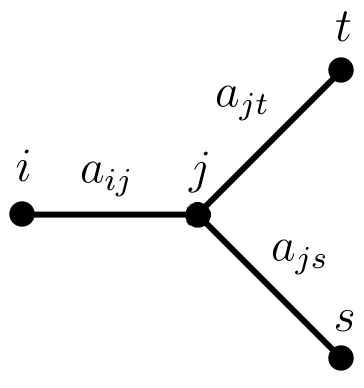}
\caption{Fragment of the matrix graph.}\label{fig}
\end{figure}

We illustrate this on the matrix graph. Consider Figure 1: for a given vertex $j$ in expression (10), there are two vertices---$t$ and $s$---that satisfy condition (9). After the first iteration, the value of the nonzero component $b_{0i}$ of the initial vector $b_0$ is propagated, according to (9), to the component $b_{1j}$. After the second iteration, by (10) and (11), this influence is further propagated to the components $b_{2t}$ and $b_{2s}$, corresponding to vertices $t$ and $s$ of the matrix graph such that $(i,t), (i,s)\!\in\!\mathcal{N}_2$.

Thus, within two iterations, an update to a single component of $b_0$---intro\-duced during the multiplication $A\cdot b_0$---propagates to those components of the vector $b_2$ that correspond to vertices at graph distance $2$ from the vertex associated with the updated component of $b_0$, i.e., the component of $b_1$ resulting from that update. This holds for every component of $b_1$ that was updated in the previous iteration.

When computing the product $ A\cdot b_1$, multiplication operations of the form $a_{ij}\cdot~b_{1j}$ are again performed for all edges $(i,j)\!\in \!\mathcal{N}_1$. This includes the components $b_{1j}$ that were updated during the first iteration, and---just as in the first iteration---each such multiplication is accompanied by an addition, potentially introducing two rounding errors per nonzero entry. 

In addition, due to the propagation of information over walks of length two, new interactions corresponding to vertex pairs at distance two---i.e., $(i,t)\!\in\!\mathcal{N}_2$---also entail arithmetic operations, yielding an additional rounding errors.

Thus, the influence of the initial perturbations---and any associated rounding errors---spreads outward in the graph, reaching vertices at distance $k$ after $k$ iterations. Similarly, at the $k$-th iteration of process (8), for each nonzero component $b_{0i}$ of the initial vector $b_0$, the components $b_{kj}$ of the vector $b_k$ corresponding to vertices $j$ such that there exists walk of the lenght $k$ between them will be updated with rounding error that propagates through this walk.  

Moreover, as noted above, all components of the vector that were modified in previous iterations continue to be updated at the $k$-th iteration as well. Thus, the $k$-th iteration involves both the propagation of information over new walks of length $k$ (captured by $\mathcal{N}_k$) and the repeated updating of components already affected in earlier steps through shorter walks.

Consequently, all else being equal---or when other relevant parameters of the matrices are similar---matrices whose graphs have a small diameter reach the maximal level of rounding error accumulation much faster than matrices with a large graph diameter for the iterative process (8).\\

\subsection{Propagation of errors during CG iterative process}

A single CG iteration involves the following operations:

\begin{itemize}
    \item[(a)] computation of the matrix-vector product $ A\cdot d_k $, where $ d_0\!=\!r_0 = b\!-\!A x_0 $;
    \item[(b)] evaluation of three inner products---$r_k^{\top}r_k$, $d_k^{\top}A d_k$, and $r_{k+1}^{\top}r_{k+1}$---used to compute the coefficients $\alpha_k$ and $\beta_k$;
    \item[(c)] three scalar-vector multiplications and three vector additions to update $x_k$, $r_k$, and $d_k$.
\end{itemize}

During the computation of the products $A d_k$ in the CG iterative process, rounding errors propagate across the components of the vectors $d_k$ in exactly the same manner as in the iterative process (8). In addition to these errors, the operations listed in items (b) and (c) introduce their own rounding errors.

When computing inner products, rounding errors accumulated in different components of the approximate solution vector, the residual, and the search direction are summed together and thus ``mixed,'' producing {\it global} errors in scalar coefficients $\alpha_k$ and $\beta_{k+1}$  that are then applied uniformly to all components during vector updates. However, the local structure of error propagation---i.e., which components influence which others and at what iteration---is governed by the topology of the matrix graph. Consequently, the graph induces deviations from the uniformity imposed by these global scalar operations.

Thus, during the execution of CG, rounding errors arise from two qualitatively distinct sources. First, {\it local errors} occur in the matrix-vector product $A\cdot d_k$. These errors propagate along the connectivity structure of the matrix graph and depend on its topology---particularly on the graph diameter. They determine how rapidly information (and numerical perturbations) spread among vector components. Second, {\it global errors} are introduced by inner products and vector updates. Since these steps employ common scalar coefficients $\alpha_k$ and $\beta_{k+1}$, the resulting errors affect all vector components uniformly. These global errors do not depend on the sparsity pattern of $A$ and contribute an additional, homogeneous amount to the total rounding error at each iteration.

Consequently, the overall accumulation of rounding errors can be viewed as the sum of a topology-dependent (local) component and a global component  independent of matrix sparsity pattern. The topology-dependent component is determined by the matrix structure, and in particular by its diameter.\\

\section{Computing Pseudo-Diameter of the Graph Matrix}

Known algorithms for computing the exact diameter $\ell$ of a sparse graph have computational complexity $O(mn)$. Consequently, the diameter of the matrix graph is not a quickly computable parameter. Instead of computing the exact value $\ell$, we estimate it by performing two breadth-first searches (BFS) on the graph $G$ and use the result as the {\it pseudo-diameter} $\tilde{\ell}$. The computational complexity of BFS on a sparse graph is $O(n)$ \cite{Burkhardt}. The 2BFS algorithm for approximating the pseudo-diameter $\tilde{\ell}$ of the matrix graph $G$ is as follows.

\begin{codebox}
\Procname{$\proc{\bf 2BFS}$}
\li \mbox{\bf Input:}\ $G$.\\

\li Choose an arbitrary vertex $s$ as the start vertex.\\

\li Traverse $G$ using BFS to find a vertex $t$ farthest from $s$.\\

\li Using $t$ as the start vertex, traverse $G$ with BFS to find a vertex $u$\\ farthest from $t$.\\

\li The distance $\tilde{\ell}\!=\! l(t,u)$ is the pseudo-diameter of the graph $G$.\\ 

\li \mbox{\bf Output:}\ $\tilde{\ell}$.

\end{codebox}

In our experiments --- whose results are presented below and in the Appendix --- we performed the 2BFS algorithm with a single choice of the initial vertex $t$ to estimate the pseudo-diameter. The accuracy of the diameter estimate $\tilde{\ell}$ can be improved by running 2BFS multiple times from different starting vertices and selecting the largest pseudo-diameter $\tilde{\ell}$ obtained across all runs.

\bigskip

\noindent{\bf Remark 1.} The pseudo-diameter of a graph, especially when computed via multiple repeated runs of the 2BFS algorithm, can be viewed as an approximation of the effective diameter---the distance that separates at least a fixed proportion of vertex pairs (the threshold value defining this proportion is typically set to about $0.9$ of the total number of vertex pairs). For matrices whose graphs have a less regular structure than the extended-star graphs (one of the graph classes considered in this work) or the graphs of Laplacians arising from partial differential equations, the effective diameter proves to be a more reliable indicator of the number $N_1$ of CG iterations that can be performed in the first stage before rounding errors become too large.

\paragraph{Block matrices and matrices with disconnected graphs.} The sparse banded matrices used in our experiments may have disconnected graphs. In general, if the graph of a matrix is disconnected, then the matrix is either block-diagonal or can be transformed into block-diagonal form by a simultaneous permutation of its rows and columns (i.e., by a symmetric permutation).

Let $G_1,\ldots,G_K$ be the connected components of the graph $G$ associated with the matrix $A$, and let $\ell_i$ denote the diameter of $G_i$, $i\!=\!\overline{1,K}$. After $\ell_i$ iterations, the submatrix of $A$ corresponding to $G_i$ becomes dense (i.e., all its entries are typically nonzero).

For a block-diagonal matrix, let $\ell_{\max}\!=\!\max\{\ell_i\}$. Then, after $\ell_{\max}$ iterations, all diagonal blocks are expected to be filled with nonzero entries (with high probability, under generic non-cancellation assumptions). For matrices with disconnected graphs, we define the pseudo-diameter parameter as $\tilde{\ell}_{\max}\!=\!\max\{\tilde{\ell}_i\}$, where $\tilde{\ell}_i$ is the pseudo-diameter of the connected component $G_i$, computed via the 2BFS procedure.

\section{The Optimal Tolerance $\hat{\varepsilon}_1$ in the First Stage\\ of the {\bf Algorithm~I} and the Matrix Graph Diameter}
With limited floating-point precision (i.e., a short significand), rounding errors accumulate more rapidly during CG iterations. Consequently, if a set of quickly computable parameters is sufficient to predict the optimal tolerance $\hat{\varepsilon}_1$ (4), then the following condition must hold: the less favorable the parameter values for a given matrix, the smaller the number $N_1$ of CG iterations that can be performed in the first stage before rounding errors become too large. In other words, for matrices with ``unfavorable'' parameter values, only a limited number of first-stage iterations can be executed while still achieving a tolerance sufficient to minimize the number $N_2$ of second-stage iterations.

At the same time, for matrices with ``favorable'' parameter values, rounding errors accumulate more slowly, allowing more iterations to be performed in the first stage and thus enabling the use of a tighter switching tolerance $\varepsilon_1$.


\begin{figure}[htbp]
\begin{center}
\includegraphics[width=100mm]{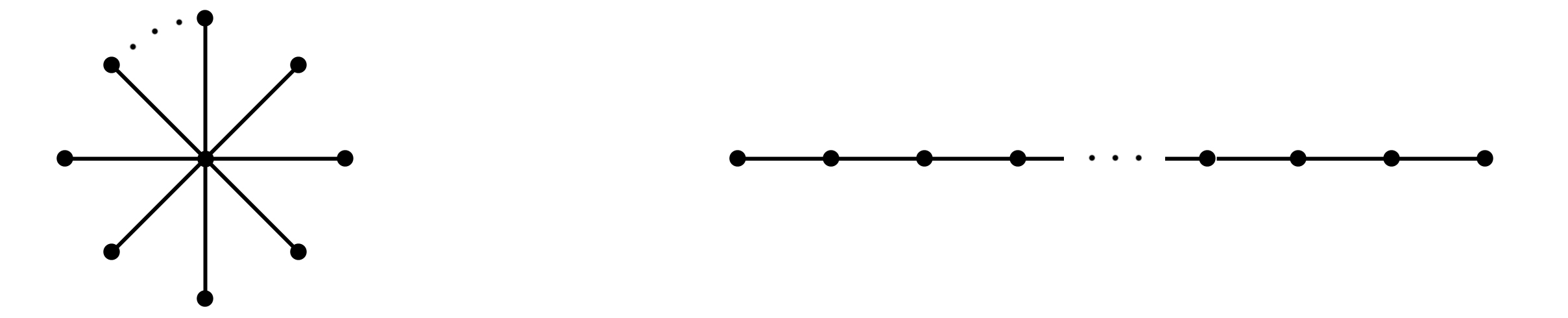}
\caption{Star and path graphs.}\label{fig1}
\end{center}
\end{figure}

Note that classifying matrix parameter values as ``favorable'' or ``unfavorable'' is not decisive for achieving the lowest total computational complexity across both stages of the algorithm. Indeed, for certain matrices, only a small number of iterations may be required in total to reach the tightest achievable residual tolerance in single precision --- rendering the distinction between ``favorable'' or ``unfavorable'' less relevant in such cases. Examples include modified adjacency matrices of star graphs with strong diagonal dominance.

In this work, we use the terms ``favorable'' and ``unfavorable'' specifically with respect to the rate of rounding error accumulation: parameters are deemed ``favorable'' if they lead to slower error growth, and ``unfavorable'' if they cause rapid error amplification.


\begin{figure}[H]
\centering
\includegraphics[width=125mm]{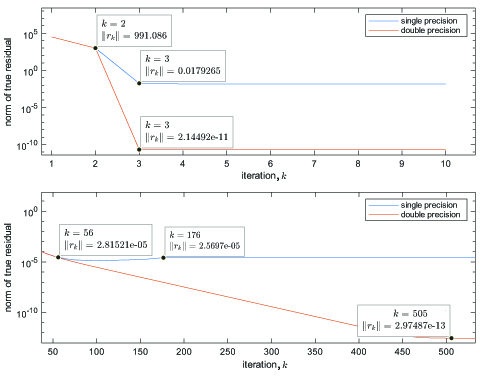}
\caption{True residual norms in CG computed in single and double precision.}\label{fig2}
\end{figure}

Consider examples of matrices whose graphs attain the minimal and maximal possible diameters: the star and the path on $n\!=\!1001$ vertices (see Fig.~\ref{fig1}). The adjacency matrices of these graphs were modified to ensure diagonal dominance by adjusting their diagonal entries as follows: $$a_{ii} = 1.001 \cdot \sum_{j \ne i} a_{ij}.$$

The first of these graphs has the minimal possible diameter $\ell\!=\!2$, and already the second power of the matrix $A$ becomes fully dense (all entries nonzero). The second graph has the maximal possible diameter $\ell\!=\!n\!-\!1$, and only the $(n\!-\!1)$-st power of the matrix $A$ becomes fully dense.

Figure 2 shows the stagnation of CG convergence for matrices whose graphs are the star and the path with $n\!=\!1001$, for which $\ell\!=\!2$ and $\ell\!=\!1000$, respectively.  

In the first case, rounding errors accumulate so rapidly that by the 2nd iteration the achievable accuracy---measured in the norm of the true residual (i.e., the absolute residual $r_k\!=\!b\!-\!A x_k$, not the recursively computed residual used internally by CG)---differs by nine orders of magnitude between single and double precision. No further improvement in accuracy occurs after $k\!=\!3$.  

In contrast, for the matrix with the path graph, the solutions in single and double precision remain nearly indistinguishable until $k\!=\!56$, and stagnation of convergence occurs only at $k\!=\!176$ for single precision and at $k\!=\!506$ for double precision.

Thus, we can perform significantly more iterations using lower precision (single precision) for the path graph and still achieve higher computational accuracy---i.e., a smaller true residual. Likewise, even in double precision, the path-graph matrix allows us to attain a solution whose true residual is about two orders of magnitude smaller than that obtained for the star-graph matrix.

\bigskip


\begin{figure}[htbp]
\centering
\includegraphics[width=50mm]{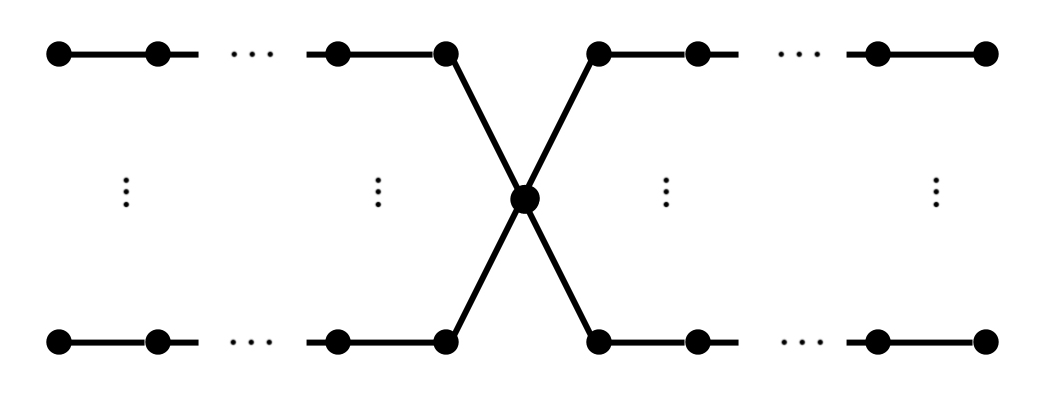}
\caption{Extended star graph.}\label{fig4}
\end{figure}

In an extended star graph, a central node is connected to several paths (\textit{rays}) of equal length (Fig.~\ref{fig4}). An examination of extended star matrices with $n\!=\!1001$ vertices shows that the values of $N_1$ range from $4$ to $31$, depending on the graph diameter.

We adopt the following classification of matrices based on their condition number $\kappa$:
\begin{itemize}
\item[1)] {\it well-conditioned} if $\kappa\!\le\!10$;
\item[2)] {\it satisfactorily conditioned} if $10\!<\!\kappa\!\le\!100$;
\item[3)] {\it moderately ill-conditioned} if $100\!<\!\kappa\!\le\!1000$;
\item[4)] {\it ill-conditioned} if $\kappa\!>\!1000$.
\end{itemize}

We evaluate the computational efficiency of {\bf Algorithm I} for solving a linear system with matrix $A$ using the metric $$E(A)=\biggl (1-\frac{I(A)}{I_{\mathrm{dp}}(A)}\biggr )\times 100\%,\eqno(12)$$ where $I(A)$ is the number of iterations computed according to (4) with $\omega\!=\!1/3$ when executing {\bf Algorithm I} with the~optimal value of $\varepsilon_1$: $\varepsilon_1\!=\hat{\varepsilon}_1\!$, and $I_{\mathrm{dp}}(A)$ denotes the number of iterations required by CG to solve the same linear system with the same stopping tolerance $\varepsilon_2$ entirely in double precision. The~value $E(A)$ represents the percentage reduction in iteration count (and thus in computational cost, under the cost model (4)) achieved by {\bf Algorithm~I} compared to double-precision CG.

Tables 1--3 report the efficiency of the algorithm with the optimal choice of~$\varepsilon_1$ for the matrix samples used in the experiments.

Each table is organized into three column groups. The first group contains the following aggregate statistics for matrices whose condition number falls within the specified range:

\begin{itemize}
\item[--] $\mbox{\sffamily M}E$: the average efficiency, computed as the arithmetic mean of $E(A)$ over all matrices $A$ in the sample with $\kappa$ in the given range;
\item[--] $\mbox{\sffamily M}\tilde{\ell}$: the average pseudo-diameter of these matrices;
\item[--] $N$: the total number of such matrices in the sample. 
\end{itemize}

The second and third groups each contain two columns, corresponding to subsets of the above sample: 

\begin{itemize}
\item[--] Second group: $\mbox{\sffamily M}E$ and $N$ for matrices for which the effficiency is worse than the average for matrices from the sample with a condition number from the given range;
\item[--] Third group: values of $\mbox{\sffamily M}\tilde{\ell}$ and $N$ for which the efficiency is better than the average for matrices from the sample with a condition number from the given range. 
\end{itemize}

\begin{table}[H]
\begin{center}
\caption{Matrices with extended star graphs with added random edges. Pseudo-diameter of the matrix graph and the efficiency of {\bf Algorithm~I} with the optimal value of $\varepsilon_1$. Sample size is $6400$.}

\bigskip

\begin{tabular}{|c||c|c|c||c|c||c|c|}
\hline  
$\kappa$ &\quad $\mbox{\sffamily M}E$\,\qquad & $\mbox{\sffamily M}\tilde{\ell}$ & $N$ & $\mbox{\sffamily M}\tilde{\ell}$ & $N$ & $\mbox{\sffamily M}\tilde{\ell}$ & $N$\\  
\hline
$\kappa\!\le\!10$          & 11.31\%  & 124.8 & 994 & 111.9 & 351  & 131.9 & 643\\
\hline
$10\!<\!\kappa\!\le\!100$  & 22.21\%  & 96.2  & 2206 & 70.0 & 1006 & 118.2 & 1200\\  
\hline
$100\!<\!\kappa\!\le\!1000$& 21.36\%  & 29.7  & 2325 & 11.8 & 1024 & 43.8  & 1301\\  
\hline
$\kappa\!>\!1000$          & 8.69\%   & 4.6   & 875  & 2.6  & 254  & 5.4   & 621\\  
\hline  
\end{tabular}  
\end{center}  
\end{table}

\begin{table}[H]
\begin{center}
\caption{Matrices with random sparse graphs. Impact of the pseudo-diameter of the matrix graph and the efficiency of {\bf Algorithm~I} with the optimal value $\varepsilon_1$. Sample size is $4800$.}

\bigskip

\begin{tabular}{|c||c|c|c||c|c||c|c|}
\hline  
$\kappa$ &\quad $\mbox{\sffamily M}E$\,\qquad & $\mbox{\sffamily M}\tilde{\ell}$ & $N$ & $\mbox{\sffamily M}\tilde{\ell}$ & $N$ & $\mbox{\sffamily M}\tilde{\ell}$ & $N$\\  
\hline
$\kappa\!\le\!10$          & 27.46\% & 3.9  & 228  & 3.9  & 118  & 4.0 & 110\\
\hline
$10\!<\!\kappa\!\le\!100$  & 24.92\% & 10.3 & 1799 & 7.8  & 823  & 12.4 & 976\\  
\hline
$100\!<\!\kappa\!\le\!1000$& 22.58\% & 42.7 & 2594 & 36.9 & 1232 & 47.8 & 1362\\  
\hline
$\kappa\!>\!1000$          & 15.19\% & 62.9 & 179  & 46.0 & 89   & 79.6 & 90\\  
\hline  
\end{tabular}  
\end{center}  
\end{table}

\begin{table}[H]
\begin{center}
\caption{Sparse banded matrices. Impact of the pseudo-diameter of the matrix graph and the efficiency of {\bf Algorithm~I} with the optimal value of $\varepsilon_1$. Sample size is $6000$.}

\bigskip

\begin{tabular}{|c||c|c|c||c|c||c|c|}
\hline  
$\kappa$ &\quad $\mbox{\sffamily M}E$\,\qquad & $\mbox{\sffamily M}\tilde{\ell}$ & $N$ & $\mbox{\sffamily M}\tilde{\ell}$ & $N$ & $\mbox{\sffamily M}\tilde{\ell}$ & $N$\\  
\hline
$\kappa\!\le\!10$          & 19.42\% & 52.9 & 5333 & 66.5  & 2610 & 39.7 & 2723\\
\hline
$10\!<\!\kappa\!\le\!100$  & 18.16\% & 256.7 & 459 & 271.1 & 266  & 237.0 & 193\\  
\hline
$100\!<\!\kappa\!\le\!1000$& 15.85\% & 225.0 & 135 & 213.4 & 96   & 253.6 & 39 \\  
\hline
$\kappa\!>\!1000$          & 13.36\% & 67.0  & 73  & 64.1  & 49   & 72.8  & 24\\  
\hline  
\end{tabular}  
\end{center}  
\end{table}

The results presented in Tables 1 and 2 show that for the samples for which they were obtained, the following statement is true:

\bigskip

{\it For satisfactorily conditioned, moderately ill-conditioned, and ill-conditioned matrices, the efficiency of {\bf Algorithm~I} is primarily determined by the diameter of the matrix graph.}    

\bigskip

This statement can also be formulated as follows: the efficiency depends on the graph diameter. This dependence is also observed for moderately ill-conditioned and ill-conditioned banded matrices (Table 3), but not for well-conditioned sparse banded matrices. A well-conditioned matrix has a beneficial effect: the accumulation of rounding errors proceeds more slowly. The sample of sparse banded matrices we considered is generally well-conditioned: approximately $97\%$ are either well- or satisfactorily conditioned, with more than $81\%$ falling into the well-conditioned category.

\bigskip

\noindent{\bf Remark 2.} For matrices with large condition numbers and/or small graph diameters, cases may arise in which no iteration savings are achieved, according to criterion (12). That is, $I(A)\!>\!I_{\mathrm{dp}}(A)$, which corresponds to a negative efficiency $E(A)\!<\!0\%$. This happens because the approximate solutions of the systems with these matrices, obtained at the first stage, deviate significantly from the solutions with the required stopping tolerance. Consequently, in the second stage, more number of iterations are required to achieve it.\\

\subsection{Propagation of Rounding Errors When Precondition is Applied}

Preconditioning is used when solving linear systems with ill-conditioned matrices. Let~$M$ be the preconditioning matrix, whose inverse is applied to the residual vectors during CG iterations. The preconditioned CG algorithm is as follows. 

\begin{codebox}
\Procname{$\proc{\bf CG with preconditioning}$}
\li \mbox{\bf Input:}\ $A$, $b$, $x_0$, $\varepsilon$.\\

\li $i\leftarrow 0$, $r_0\leftarrow M^{-1}(b-Ax_0)$, $d_0\leftarrow r_0$.\\

\li\While $\|Ax_i-b\|>\varepsilon$:
   \Do 

\li $\alpha_i\leftarrow r_i^{\top}M^{-1}r_i/d_i^{\top}Ad_i$;

\li $x_{i+1}\leftarrow x_i+\alpha_id_i$;

\li $r_{i+1}\leftarrow r_i-\alpha_iAd_i$;

\li $\beta_{i+1}\leftarrow r_{i+1}^{\top}M^{-1}r_{i+1}/r_i^{\top}M^{-1}r_i$;

\li $d_{i+1}\leftarrow M^{-1}r_{i+1}+\beta_{i+1}d_i$;

\li $i\leftarrow i+1$;
    \End

\li $\tilde{x}\leftarrow x_i$.\\

\li \mbox{\bf Output:}\ Approximate solution $\tilde{x}$. 		

\end{codebox}

The question of how the diameter and other topological characteristics of the matrix graph, which we consider, affect the convergence of the conjugate gradient method when preconditioning is applied lies beyond the scope of this work. We only note that in the case of preconditioning, the relationship between the rate of loss of orthogonality of residuals and conjugacy of search directions and the diameter (effective diameter, pseudo-diameter) is preserved provided that the preconditioner does not disrupt the structure of information and rounding error propagation determined by the system matrix. This condition is satisfied if multiplying the residual by the inverse of the preconditioner matrix---$M^{-1}r_{i+1}$---in the recurrence for the search direction $$d_{i+1}\leftarrow M^{-1}r_{i+1}+\beta_{i+1}d_i$$ does not introduce additional information propagation beyond that occurring in the unpreconditioned CG method. Examples of such preconditioners include Jacobi preconditioning for arbitrary symmetric matrices, which is effective for diagonally dominant matrices, as well as incomplete Cholesky preconditioning (IC(0)) for moderately ill-conditioned matrices.

In our experiments that we consider further, we omitted preconditioning. CG without preconditioning requires more iterations at both stages than its preconditioned counterpart, yielding a~broader range of possible $(N_1,N_2)$ pairs of the iteration counts in the two stages. This increased variability makes the classification of matrices by $\hat{\varepsilon}_1$ more challenging for the relatively small-dimensional problems considered in this work and allows us to demonstrate the effectiveness of our approach more objectively.\\

\section{The Diameter of the Matrix Graph\\ and Rounding Error Growth\\ in Other Iterative Methods for Solving Linear Systems}

Let us examine the diameter of the matrix graph and the growth of rounding errors in some other iterative methods for solving linear systems. 

Let $T$ be the transition matrix of the iterative method: $$x_{k+1}=Tx_k+c,$$ where $x_k$ is the approximate solution obtained at the $k$-th iteration, and $c$ is a~vector that depends on the method and the right-hand side $b$ of the system. If~the transition matrix $T$ has the same graph as the system matrix, then its diameter also influences the accumulation of rounding errors during the iterations.

In particular, this holds for the Jacobi method, whose transition matrix shares the same sparsity graph as the system matrix. The transition matrix is given by $$T=D^{-1}(L+U),$$ where $D\!=\!\mbox{diag}(a_{11},\ldots, a_{nn})$ is the diagonal part of $A$, $L$ is the strictly lower triangular part (with entries $a_{ij}$ for $i\!>\!j$), and $U$ is the strictly upper triangular part (with entries $a_{ij}$ for $j\!>\!i$), so that $A\!=\!L+D+U$. The Jacobi iteration reads $$x_{k+1}=D^{-1}((L+U)x_k+b),$$ and the error propagation is described by $$e_k =(D^{-1}(L+U))^ke_0.$$ Since $T$ and $A$ have identical sparsity patterns, the diameter of the matrix graph directly influences the error propagation in the Jacobi method.

At the same time, for another simple iterative method --- the Gauss--Seidel method --- the transition matrix has a graph different from that of the system matrix: $$T=(D-L)^{-1}U,\quad e_k=((D-L)^{-1}U)^ke_0,$$ and no such dependence arises. The same holds for other relaxation methods, whose transition matrices generally take the form $$T=(D-\gamma L)^{-1}((1-\gamma)D+\gamma U),\quad \gamma\!\in\!(0,2).$$

Tables 4 and 5 report the number of iterations required by the conjugate gradient (CG) and Jacobi methods in single precision ($N_1$) and double precision ($N_2$) to obtain an approximate solution with stopping tolerance $\varepsilon_2\!=\!10^{-10}$, for $\varepsilon_1\!=\!10^{-t}$, $t\!=\!\overline{2,7}$. Here, $I(A)$ denotes the complexity estimate defined in (4).

Table 6 presents the corresponding values of $I(A)$, $N_1$, and $N_2$ for the Gauss--Seidel method. For comparison, we use linear systems whose coefficient matrices are derived from the adjacency matrices of star and path graphs. The parameter $\omega$ is set to $1/3$, $\mu\!=\!1.1$.

The dashes in Table 5 for the Jacobi method indicate that the tolerance $\varepsilon_1\!=\!10^{-t}$ cannot be achieved for the corresponding values of $t$, as rounding errors become too large in single precision. 

For both the Jacobi method and CG, the optimal number of iterations in the first stage of the algorithm --- determined by the growth of rounding errors during single-precision computations --- depends on the diameter of the system matrix graph (see the marked rows in the tables). In contrast, for the Gauss--Seidel method, whose transition matrix has a graph different from that of the system matrix, no such dependence is observed.

\begin{table}[H]
\begin{center}
\caption{Number of single- and double-precision iterations of {\bf Algorithm~I} using CG for different values of $\varepsilon_1$, with system matrices derived from star and path graphs ($n\!=\!1001$).}

\bigskip

\begin{tabular}{|c|c|c|c||c|c|c|}
\hline  
&\multicolumn{3}{|c||}{Star}&\multicolumn{3}{|c|}{Path}\\  
\cline{2-7}  
$\varepsilon$ & $I(A)$ & $N_1$ & $N_2$ & $I(A)$ & $N_1$ & $N_2$\\  
\hline  
$10^{-2}$& 4.3  & 4  & 3  & 18.3 & 7  & 16\\  
\cline{2-4}
$10^{-3}$& 4.3  & 4  & 3  & 17.0 & 9  & 14\\  
$10^{-4}$& 4.7  & 5  & 3  & 15.7 & 11 & 12\\  
\cline{5-7}
$10^{-5}$& 5.0  & 6  & 3  & 15.0 & 12 & 11\\  
\cline{5-7}
$10^{-6}$& 5.7  & 8  & 3  & 15.0 & 12 & 11\\  
$10^{-7}$& 5.7  & 8  & 3  & 15.3 & 16 & 10\\  
\hline  
\end{tabular}  
\end{center}  
\end{table}

\begin{table}[H]
\begin{center}
\caption{Number of single- and double-precision iterations of {\bf Algorithm I} using Jacobi's method for~different values of $\varepsilon_1$, with system matrices derived from star and path graphs ($n\!=\!1001$).}

\bigskip

\begin{tabular}{|c|c|c|c||c|c|c|}
\hline  
&\multicolumn{3}{|c||}{Star}&\multicolumn{3}{|c|}{Path}\\  
\cline{2-7}  
$\varepsilon$ & $I(A)$ & $N_1$ & $N_2$ & $I(A)$ & $N_1$ & $N_2$\\  
\hline  
$10^{-2}$& 34.7 & 20 & 28 & 32.3 & 16 & 27\\  
\cline{2-4}
$10^{-3}$& \underline{32.3} & 22 & 25 & 30.3 & 19 & 24\\  
\cline{2-4}
$10^{-4}$& 32.3 & 22 & 25 & 28.3 & 22 & 21\\  
\cline{5-7}
$10^{-5}$&   -- & -- & -- & \underline{7.7} & 26 & 1\\  
\cline{5-7}
$10^{-6}$&   -- & -- & -- & 9.7 & 26 & 1\\  
$10^{-7}$&   -- & -- & -- & 9.7 & 26 & 1\\  
\hline  
\end{tabular}  
\end{center}  
\end{table}

\begin{table}[H]
\begin{center}
\caption{Number of single- and double-precision iterations of {\bf Algorithm I} using the~Gauss-Seidel method for~different values of $\varepsilon_1$, with system matrices derived from star and path graphs ($n\!=\!1001$).}

\bigskip

\begin{tabular}{|c|c|c|c||c|c|c|}
\hline  
&\multicolumn{3}{|c||}{Star}&\multicolumn{3}{|c|}{Path}\\  
\cline{2-7}  
$\varepsilon$ & $I(A)$ & $N_1$ & $N_2$ & $I(A)$ & $N_1$ & $N_2$\\  
\hline  
$10^{-2}$& 17.0 & 9  & 14 & 16.7 & 8  & 14\\  
$10^{-3}$& 14.3 & 10 & 11 & 16.0 & 9  & 13\\  
$10^{-4}$& 5.7  & 11 & 2  & 14.7 & 11 & 11\\  
$10^{-5}$& 5.7  & 11 & 2  & 13.3 & 13 & 9\\  
\cline{2-7}
$10^{-6}$& 5.0  & 12 & 1  & 6.7  & 14 & 2\\  
\cline{2-7}
$10^{-7}$& 5.0  & 12 & 1  & 6.7  & 14 & 2\\  
\hline  
\end{tabular}  
\end{center}  

\end{table}

\section{The Average Rate of Residual Norm Decay\\ as a Feature for Predicting $\hat{\varepsilon}_1$}

In addition to the matrix size $n$, the number of nonzeros $m$, and the pseudo-diameter $\tilde{\ell}$ of the matrix graph, we also use the average rate $v$ of residual norm decay during the early iterations of CG in single precision to predict the optimal tolerance $\hat{\varepsilon}_1$.

The ratio of residuals at consecutive CG iterations depends on the condition number $\kappa$ of the matrix. Specifically, for the residual norms of the initial approximate solution and the solution at the $i$-th iteration, the following asymptotic bound --- analogue to the error estimate (6), but in the Euclidean norm --- holds: $$\frac{\|r_i\|_2}{\|r_0\|_2}\le 2\biggl (\frac{\sqrt{\kappa}-1}{\sqrt{\kappa}+1}\biggr )^i.$$ Consequently, the per-iteration decay rate $v_i$ satisfies $$v_i=\frac{\|r_i\|_2}{\|r_{i-1}\|_2}\le\frac{\sqrt{\kappa}-1}{\sqrt{\kappa}+1}.$$

The residual decay rate varies most significantly during the first few iterations, as CG rapidly eliminates error components associated with the extreme eigenvalues of the matrix. Moreover, the initial residual is large, so the earliest corrections yield the greatest reduction in error. By combining the average of~$v_i$ over the early single-precision iterations with the other parameters described above, we can predict $\hat{\varepsilon}_1$ with sufficient accuracy to achieve near-optimal algorithmic efficiency.

Thus, the fourth parameter used to predict $\hat{\varepsilon}_1$ for the input matrix is the average rate v of residual norm decay during the early iterations of CG in single precision: $$v=\frac{1}{k_0}\sum\limits_{i=1}^{k_0}\frac{\|r_i\|}{\|r_{i-1}\|}=\frac{1}{k_0}\sum\limits_{i=1}^{k_0}v_i.\eqno(13)$$ Here, $k_0$ is an algorithmic parameter whose value is determined in advance by sampling matrices of the type to which the algorithm is applied. 

\bigskip

\noindent{\bf Remark 3.} Determining the parameter $v$ does not require any computations beyond those already performed by the algorithm. After computing $v$ from the early single-precision iterations and selecting $\varepsilon_1$ via the $k$-nearest neighbors method, we continue iterating in single precision, using the solution obtained at iteration $k_0$ as the initial guess. The residual norms $\|r_i\|_2$ for $i\!\le\!k_0$ are computed in double precision. 

\bigskip

Thus, we use the following precision-switching algorithm. 

\begin{codebox}
\Procname{$\proc{\bf Algorithm II}$}
\li \mbox{\bf Input:}\ $A$, $b$, $x_0$, $\varepsilon_2$, $k_0$.\\

\li Convert $A$ and $b$ to single precision: $\tilde{A}\!\leftarrow\!\texttt{float32}(A)$, $\tilde{b}\!\leftarrow\!\texttt{float32}(b)$.\\

\li Set the initial approximate solution to zero: $\tilde{x}_0\!\leftarrow\!0$.\\

\li Perform $k_0$ single-precision CG iterations and determine the value of $v$.\\

\li Determine $\varepsilon_1$ by classifying A via its feature vector $\chi(A)\!=\!(n,m,\tilde{\ell},v)$.\\

\li Using CG, compute an approximate solution  $\tilde{x}$ of the system $\tilde{A}x\!=\!\tilde{b}$\\ 
with tolerance $\varepsilon_1$, using the solution obtained at iteration $k_0$\\ as the initial guess.\\

\li Convert $\tilde{x}$ to a double-precision vector $x_0$: $x_0\!\leftarrow\!\texttt{float64}(\tilde{x})$.\\

\li Using the initial approximate solution $x_0$, find an approximate solution $\breve{x}$\\ of the system $Ax\!=\!b$ with a stopping tolerance of $\varepsilon_2$.\\

\li \mbox{\bf Output:}\ Refined approximate solution $\breve{x}$.

\end{codebox}

\section{The Complexity of Computing the Parameters\\ used to Predict $\hat{\varepsilon}_1$}

Let $C(\mbox{CG}_{\mathrm{dp}})$ denote the computational complexity, measured in clock cycles, required to solve a linear system with residual tolerance $\varepsilon_2$ using CG executed entirely in double precision, and let $C(\mbox{Alg})$ denote the complexity (in clock cycles) required to solve the same linear system using {\bf Algorithm~II}. The total cost of {\bf Algorithm~II} is given by $$C(\mbox{Alg})=C(\mbox{EC}^{\mathrm{Alg}})+C(\mbox{CG}_{\mathrm{sp}}^{\mathrm{Alg}})+C(\mbox{CG}_{\mathrm{dp}}^{\mathrm{Alg}}),$$ where 
\begin{itemize}
\item[--] $C(\mbox{EC}^{\mathrm{Alg}})$ (\textit{EC is an abbreviation for Extra Computations}) is the~com\-ple\-xity required to compute the matrix parameters used to predict $\hat{\varepsilon}_1$, and for the prediction itself, 
\item[--] $C(\mbox{CG}_{\mathrm{sp}}^{\mathrm{Alg}})$ and $C(\mbox{CG}_{\mathrm{dp}}^{\mathrm{Alg}})$ is the complexity of executing CG with single and double precision in the stages of {\bf Algorithm II}.
\end{itemize}

To evaluate the efficiency of {\bf Algorithm II}, it is necessary that for a small value of $\alpha$, $$C(\mbox{EC}^{\mathrm{Alg}})\le\alpha C(\mbox{CG}_{\mathrm{dp}}).\eqno(14)$$ In the estimates obtained below, we proceed from the value $\alpha\!=\!0.01$; that is, the com\-pu\-tational complexity $C(\mbox{EC}^{\mathrm{Alg}})$, expressed by us in clock cycles, must not exceed $1\%$ of the com\-pu\-tational complexity $C(\mbox{CG}_{\mathrm{dp}})$.

As noted above, computing the values of the parameters $n$, $m$ and $v$ does not require time additional to the running time of {\bf Algorithm II}. The extra cost $C(\mbox{EC}^{\mathrm{Alg}})$ arises only from:
\begin{itemize}
\item[1)] computing the pseudo-diameter $\tilde{\ell}$ of the system matrix $A$;
\item[2)] classifying the matrix $A$ by its feature vector $\chi(A)\!=\!(n,m,\tilde{\ell},v)$ and predicting of the optimal tolerance $\hat{\varepsilon}_1$ based on the matrix class.
\end{itemize}

\smallskip

Let us estimate the computational cost $C(\mbox{EC}^{\mathrm{Alg}})$.\\

\subsection{Computational cost of CG in double precision (clock cycles)}

To justify the use of matrix parameters that are quickly computable, we estimate the computational cost $C(\mbox{CG}_{\mathrm{dp}})$ of CG in double precision, measured in processor clock cycles.

\smallskip

In estimating $C(\mbox{CG}_{\mathrm{dp}})$, as well as the costs of computing the pseudo-diameter $\tilde{\ell}$ and performing the $k$-nearest neighbors classification, we account only for clock cycles spent on arithmetic and logical operations. We exclude cycles associated with memory accesses and assignment operations, as their count depends on memory hierarchy organization and data layout --- factors determined by the specific hardware platform and runtime environment.

All complexity estimates for the algorithms presented below are derived for their sequential implementations. When parallel variants are employed, the estimates do not degrade; on the contrary, they become more favorable for the mixed-precision algorithm under consideration.

\smallskip

Since the algorithm is intended for large sparse matrices, we first estimate $C(\mbox{CG}_{\mathrm{dp}})$ for preconditioned CG.

\smallskip

We now estimate the computational cost of a single iteration, accounting only for arithmetic operations whose total cost scales linearly with $n$ (i.e., those that dominate the asymptotic complexity). Let 
\begin{itemize}
\item[--] $C_{+}$ be the computational complexity of addition in cycles, 
\item[--] $C_{\times}$ be the complexity of multiplication, 
\item[--] $C_{\mathrm{sc}}(n)$ be the complexity of the scalar product of two $n$-dimensional vectors,  
\item[--] $C_{\mathrm{mlt}}(m)$ be the complexity of the product of an $n$-dimensional vector by a~sparse $n\!\times\! n$ matrix $A$ with $m$ nonzero elements. 
\end{itemize}

The scalar product requires $n$ multiplications and $n\!-\!1$ additions, so
$$C_{\mathrm{sc}}(n)=nC_{\times}+(n-1)C_{+}\approx nC_{\times}+nC_{+}=n(C_{\times}+C_{+}).\eqno(15)$$ For sparse matrix-vector multiplication, each of the $m$ nonzeros contributes one multiplication and one addition (for accumulation), yielding $$C_{\mathrm{mlt}}(m)=mC_{\times}+mC_{+}=m(C_{\times}+C_{+}).\eqno(16)$$

We now analyze the cost of each step in one CG iteration with preconditioning:
\begin{itemize}
\item[--] {\bf Step 2:} (residual norm check $\|Ax_i-b\|_2>\varepsilon$): $C_{\mathrm{mlt}}(n)+nC_++C_{\mathrm{sc}}(n)$.\\
(matrix-vector product, scalar product for norm, and vector subtraction).
\item[--] {\bf Step 3:} $2(C_{\mathrm{mlt}}(m)+C_{\mathrm{sc}}(n))$. 
\item[--] {\bf Step 4:} $nC_{\times}+nC_+$. 
\item[--] {\bf Step 5:} $C_{\mathrm{mlt}}(m)+nC_{\times}+nC_+$. We exclude $C_{\mathrm{mlt}}(m)$ here, as this matrix-vector product is already computed in Step 3. 
\item[--] {\bf Step 6:} $2(C_{\mathrm{mlt}}(m)+C_{\mathrm{sc}}(n))$. 
\item[--] {\bf Step 7:} $C_{\mathrm{mlt}}(m)+n C_{\times}+nC_+$. 
\end{itemize}

Summing the dominant costs and substituting (15)--(16), the total cost of one preconditioned CG iteration is $C(\mbox{CG}_{\mathrm{iter}}^{\mathrm{c}})$ of one iteration of CG with preconditioning: $$C(\mbox{CG}_{\mathrm{iter}}^{\mathrm{c}})=6C_{\mathrm{mlt}}(m)+5C_{\mathrm{sc}}(n)+3nC_{\times}+4nC_+=$$
$$=6(mC_{\times}+mC_+)+5(n C_{\times}+nC_+)+3nC_{\times}+4nC_+=(6m+8n)C_{\times}+(6m+9n)C_+.$$

Modern processors typically require $3$--$4$ clock cycles to perform an addition ($C_{+}$) and $3$--$5$ cycles for a multiplication ($C_{\times}$), for both single- and double-precision floating-point numbers. For simplicity, we set $C_{+}\!=\!4$ and $C_{\times}\!=\!4$. Substituting into the per-iteration cost yields $$C(\mbox{CG}_{\mathrm{iter}}^{\mathrm{c}})=(6m+8n)4+(6m+9n)4=4(12m+17n).\eqno(17)$$ The total computational cost of preconditioned CG, which performs $k_{\mathrm{CG}^{\mathrm{c}}}$ double-precision iterations, is therefore $$C(\mbox{CG}_{\mathrm{dp}}^{\mathrm{c}})=C(\mbox{CG}_{\mathrm{iter}}^{\mathrm{c}})\cdot k_{\mathrm{CG}^{\mathrm{c}}}=4(12m+17n)\cdot k_{\mathrm{CG}^{\mathrm{c}}}.\eqno(18)$$

\smallskip

For CG without preconditioning, the iteration cost excludes the four sparse matrix-vector multiplications (each of cost $C_{\mathrm{mlt}}(m)$) that appear in Steps 3, 6, and 7 of the preconditioned variant. Consequently, the computational cost of unpreconditioned CG in double precision is $$C(\mbox{CG}_{\mathrm{dp}})=4(6m+17n)\cdot k_{\mathrm{CG}}\eqno(19)$$ where $k_{\mathrm{CG}}$ denotes the total number of iterations performed.\\

\subsection{Low computational cost of the pseudo-diameter\\ of a matrix graph}

To characterize matrix parameters as quickly computable and to estimate the sample size admissible for satisfying condition (14) in the $k$-nearest neighbors method, we also need non-asymptotic complexity estimates of the algorithms employed by {\bf Algorithm~II}, expressed explicitly in CPU clock cycles.

\smallskip

The optimal sequential version of algebraic BFS for sparse graphs requires $m$ elementary operations \cite{Burkhardt} --- namely logical operations and integer additions --- each of which executes in one clock cycle. Thus, the cost of the 2BFS algorithm for computing the pseudo-diameter is $C(\mbox{2BFS})\!=\!2m$.

The ratio of this cost to the runtime of unpreconditioned CG in double precision, $C(\mbox{CG}_{\mathrm{dp}}\!=\!(6c+17)n\cdot k_{\mathrm{CG}}$ (see (19)), is $$\frac{C(\mbox{2BFS})}{C(\mbox{CG}_{\mathrm{dp}})}=\frac{2m}{4(6m+17n)\cdot k_{\mathrm{CG}}}=\frac{cn}{2(6c+17)n\cdot k_{\mathrm{CG}}}=\frac{c}{2(6c+17)\cdot k_{\mathrm{CG}}},\eqno(20)$$ where $c\!=\!m/n$ denotes the average number of nonzeros per row.

For a sparse matrix with $m\!\approx\! n$, that is, for $c\!\approx\! 1$, we get $$\frac{C(\mbox{2BFS})}{C(\mbox{CG}_{\mathrm{dp}})}=\frac{1}{46\cdot k_{\mathrm{CG}}}.$$ Thus, as soon as the number of CG iterations satisfies $k_{\mathrm{CG}}\!\ge\!3$, we have $$C(\mbox{2BFS})\!<\!0.01\cdot C(\mbox{CG}_{\mathrm{dp}}),$$ i.e., the cost of computing the pseudo-diameter $\tilde{\ell}$ is less than $1\%$ of the total CG runtime. For values of $c$ such that the $1000\!\times\!1000$ matrix contains at most $10\%$ nonzeros --- i.e., for $c\!\le\!100$ --- the required number of iterations does not exceed~$8$.

Since each iteration of preconditioned CG is more expensive than that of~unpreconditioned CG, the ratio $C(\mbox{CG}_{\mathrm{dp}}^{\mathrm{c}})$ is even smaller. Consequently, for preconditioned CG, the $1\%$ threshold is reached with an even smaller number of iterations $k_{\mathrm{CG}^c}$.

Therefore, the pseudo-diameter $\tilde{\ell}$ is indeed quickly computable: for high-dimensional sparse matrices, both unpreconditioned and preconditioned CG typically require far more than the minimal number of iterations identified above, ensuring that the cost of computing $\tilde{\ell}$ remains well below $1\%$ of the total solver cost.\\

\section{Finding the Optimal Tolerance Parameter $\hat{\varepsilon}_1$\\ in the First Stage of the Algorithm\\ via a Classification Approach}

\subsection{Classification problem}

The functional dependence of the optimal tolerance $\hat{\varepsilon}_1$ --- achieved in the first stage of the algorithm --- on the matrix $A$ (or on any of its parameters) is highly complex and cannot be easily expressed explicitly as a function $\hat{\varepsilon}_1\!=\!\varphi(A)$. In particular, it is impractical to represent it explicitly as a function $\hat{\varepsilon}_1\!=\!\varphi(\chi(A))$ of the feature vector $\chi(A)=\bigl (n,m,\tilde{\ell}(A),v\bigr)$ used in our approach. This remains challenging even when attempting to learn such a mapping via neural network training, as will be discussed in more detail below.

Instead, we employ the $k$-nearest neighbors method (hereafter \emph{kNN}) to detremine $\varphi(\chi(A))$ based on the values of $\varphi$ at feature vectors $\chi$ that are close to $\chi(A)$. All feature vectors $\chi$ are componentwise normalized using the training sample. Specifically, we apply minimax normalization to map each component into the interval $(0,1]$. The same normalization --- derived from the kNN training sample --- is applied to the feature vector $\chi(A)$ of the input matrix $A$. This normalization step incurs computational cost significantly lower than $\alpha C(\mbox{CG}_{\mathrm{dp}})$ with $\alpha\!=\!0.01$ (see~(14)).

By classifying the feature vector $\chi(A)$, we effectively perform a regression of the discrete-valued function $\hat{\varepsilon}_1(A)$. To this end, we define the classes $\mathcal{C}_t$ as $$\mathcal{C}_t=\{A\in\mathbb{R}^{n\times n}\, |\, \hat{\varepsilon}_1(A)=\varepsilon_t\},\eqno(23)$$ 
where $\varepsilon_1\!=\!10^{-t}$, $t\!=\!\overline{2,7}$.

To construct a training sample of feature vectors for matrices of a given type, we first generate a set of such matrices. For each matrix $A$, we evaluate all admissible values $\varepsilon_t$ by running CG method in single precision until the residual norm falls below $\varepsilon_1$, recording the required number of iterations $N_1$. Subsequently, we execute CG in double precision to achieve the target stopping tolerance $\varepsilon_2\!=\!10^{-10}$, storing the corresponding iteration count $N_2$.

Using the complexity estimate (4), we select the optimal $\hat{\varepsilon}_1$ and thereby assign $A$ to its class $C_t$. This yields the labeled sample $$\mathcal{S}^M_0=\bigl\{\bigl (A^{(i)},\mathcal{C}(A^{(i)})\bigr )\bigr\}_{i=1}^N,$$ where $\mathcal{C}(A^{(i)})$ denotes the class (23) --- i.e., the class label --- corresponding to the matrix $A^{(i)}$. From this sample, we derive a labeled set $\mathcal{S}_0$ of normalized feature vectors: $$\mathcal{S}_0=\bigl\{\bigl (\chi^{(i)},\mathcal{C}(A^{(i)})\bigr )\bigr\}_{i=1}^N,$$ where each feature vector $$\chi^{(i)}\!=\!\bigl (n(A^{(i)}),m(A^{(i)}),\tilde{\ell}(A^{(i)}),v(A^{(i)}))\bigr )$$ inherits the class label $\mathcal{C}(A^{(i)})$ of its corresponding matrix $A^{(i)}$.\\

\subsection{Classification using kNN}

As one of the simplest machine learning algorithms, kNN is effective for a wide range of practical problems, particularly those involving small training samples. As we will demonstrate --- based on both theoretical estimates and computational experiments --- only a small sample size is required to accurately classify a matrix according to its optimal tolerance $\hat{\varepsilon}_1$.

\smallskip

Let us briefly describe the kNN. Suppose a distance metric $\rho$ is given, which quantifies the dissimilarity between objects represented by feature vectors $\chi$. Let $\mathcal{S}\!=\!\{(\psi_1,y_1),\ldots,(\psi_{N_{\mathcal{S}}},y_{N_{\mathcal{S}}})\}$ be the labeled training sample, where $y_i$ denotes the class label of $\psi_i$. For a query vector $\chi$, let $\mathcal{N}$ denote the set of its $k$ nearest neighbors in $\mathcal{S}$ under the metric $\rho$. The weighted kNN algorithm proceeds as~follows.

\begin{codebox}
\Procname{$\proc{\bf kNN}$}

\li \mbox{\bf Input:}\ $\chi$, $\mathcal{S}$.\\

\li Let the set $\mathcal{N}$ consist of the first $k$ feature vectors from the sample $\mathcal{S}$:\\
 $\mathcal{N}\leftarrow\{\psi_1,\ldots,\psi_k\}$.\\

\li \For $i\leftarrow 1$ \To $N_{\mathcal{S}}$: \Do

\li       Compute $\rho(\chi,\psi_i)$.

\li 			Find $\psi'\!\in\!\mathcal{N}$ such that $\rho(\chi,\psi')=\max\limits_{\psi\in\mathcal{N}}\rho(\chi,\psi)$.

\li       \If $\rho(\chi,\psi_i)<\rho(\chi,\psi')$:\Then

\li       $\mathcal{N}\leftarrow\mathcal{N}\setminus\{\psi'\}$;

\li       $\mathcal{N}\leftarrow\mathcal{N}\cup\{\psi_i\}$.

         \End
         
		\End\\

\smallskip
		
\li Compute the total weight of elements of each class in $\mathcal{N}$,\\ using the metric $\rho(\,\cdot\, ,\,\cdot\,)$ to assign weights.\\

\li Classify $\chi$ as belonging to the class $\mathcal{C}_{i_0}$ with the largest total weight\\ among those represented in $\mathcal{N}$.\\

\li \mbox{\bf Output:}\ Class label $i_0$.

\end{codebox}

\subsection{Localization of feature vectors from different classes}

kNN can be effectively applied to a given classification problem if the feature vectors $\chi$ are well {\it localized}, that is, the following condition holds:

\bigskip

{\it If most of the vectors near $\chi$ belong to class $\mathcal{C}_t$, then $\chi$ also belongs to $\mathcal{C}_t$.}

\bigskip


\begin{figure}[htbp]
\begin{center}
\includegraphics[width=\textwidth, height=3.5cm, keepaspectratio=false]{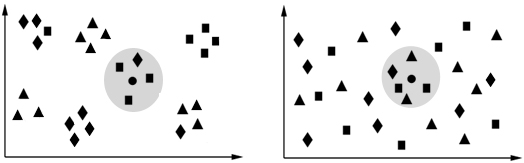}
\caption{Good and bad localization of feature vectors from the sample.}\label{fig5}
\end{center}
\end{figure}

\noindent Figure~\ref{fig5} schematically illustrates well-localized (left) and poorly localized (right) feature vectors. The elements of the three classes are represented by triangles, squares, and diamonds, respectively, while the vector to be classified is shown as a circle. In the first case, vectors from the same class tend to form compact clusters, enabling kNN to classify efficiently. The $k$ nearest neighbors --- used for classification --- are circled in gray. In the second case, such clustering is absent, which hinders effective classification with kNN.

In the experiments --- whose results are presented in the Appendix --- we quantify the {\it localization} of a sample element with respect to its class as the percentage of neighbors in $\mathcal{N}$ (the set of $k$ nearest neighbors) that belong to the~same class as the query vector. We denote this measure by $\mathrm{Loc}$. As the experimental results show, all matrix samples considered exhibit strong localization of~feature vectors across classes.\\

\subsection{Class assignment for a feature vector}

Let $\mathcal{N}\!=\!\{\psi_1,\ldots,\psi_k\}$ be the set of $k$ nearest neighbors of the classified feature vector $\chi$. That is, $\mathcal{N}$ is consists of the feature vectors from the sample $\mathcal{S}$ whose Euclidean distances to $\chi$ are among the $k$ smallest over $\mathcal{S}$.

The standard kNN (without distance weighting) assigns a class to $\chi$ based on the majority class in $\mathcal{N}$, ignoring the actual distances to $\chi$. In contrast, distance-weighted kNN takes these distances into account.

Specifically, unweighted kNN determines the class index $i_0$ as $$i_0=arg \max\limits_t N_t^{(0)},$$ where $$N_t^{(0)}=\sum\limits_{j=1}^{k}\mathbb{I}_t(\psi_j),\quad N_t^{(0)}\!\in\!\mathbb{N},\eqno(22)$$ and the indicator function is defined by $$\mathbb{I}_t(\psi_j)=\left\{
\begin{array}{ll}
1,& \psi_j\in\mathcal{C}_t,\\
0,& \psi_j\not\in\mathcal{C}_t.\\
\end{array}
\right.$$ 

We modify (22) by replacing the raw count $N_t^{(0)}$ with a distance-weighted sum: $$N_i=\sum\limits_{j=1}^k\frac{1}{\rho(\chi,\psi_j)}\mathbb{I}_t(\psi_j),\quad N_t\!\in\!\mathbb{R}_+,\eqno(23)$$ where $\rho(\chi,\psi_j)$ denotes the Euclidean distance between $\chi$ and $\psi_j$. The class $\mathcal{C}_{i_0}$ is then selected as $$i_0=arg \max\limits_t N_t.$$ Thus, neighbors farther from $\chi$ contribute less to the score $N_t$ of their respective class.

For matrix samples, distance-weighted kNN generally yields higher classification accuracy than the unweighted variant. This is confirmed by the tables ``Localization of feature vectors'' and ``kNN Accuracy'' in the Appendix. Across all experimental samples, the average accuracy $\mbox{\sffamily M}\mbox{Acc}$ of the weighted kNN consistently exceeds the average localization measure $\mbox{\sffamily M}\mbox{Loc}$.

\smallskip

The example presented in the Appendix is particularly illustrative in this regard. By enriching the base sample --- used to draw training and test sets --- with additional matrices (without increasing the size of the training set itself), we enhance its representativeness and thereby significantly improve the accuracy of kNN. Although the localization measure does not increase substantially, the feature vectors from these additional matrices are often closer to the query vector, enabling more accurate classification.

Specifically for this (see Appendix, Tables~15--24), the sample from which training feature vectors are randomly selected is augmented with vectors derived from newly generated random matrices sharing the same structural parameters. The training set size remains unchanged. Nevertheless, despite only a marginal improvement in localization, the classification accuracy of kNN increases by more than $5\%$, thereby enhancing the overall efficiency of {\bf Algorithm~II}.\\

\subsection{Computational cost of kNN and training set size}

To satisfy condition (14), the computational cost of the kNN algorithm within {\bf Algorithm~II} must not exceed a fraction $\alpha$ of  $C(\mbox{CG}_{\mathrm{dp}})$. The cost of kNN is primarily determined by two parameters: the size $N_{\mathcal{S}}$ of the training set --- used to classify the system matrix --- and the number of neighbors $k$. We denote the computational cost of kNN as $C(\mbox{kNN})\!=\!C(\mbox{kNN},N_{\mathcal{S}},k)$. For fixed matrix size $n$ and neighbor count $k$, we estimate the maximum training sample size $N_{\mathcal{S}}$ such that (14) holds.

\smallskip

We express $C(\mbox{kNN})$ in terms of clock cycles required for arithmetic and logical operations during execution: $$C(\mbox{kNN})=N_{\mathcal{S}}\cdot C(\mbox{kNN}_{\mathrm{iter}}),$$ where $C(\mbox{kNN}_{\mathrm{iter}})$ denotes the cost of a single kNN iteration.

\bigskip

During the $i$-th iteration of kNN, the following steps are performed:

\smallskip

\begin{itemize}
\item[1)] the distance $\rho(\chi,\psi_i)$ from the query vector $\chi$ to the feature vector $\psi_i$ of the training sample is computed;
\item[2)] a decision is made whether to include $\psi_i$ in $\mathcal{N}$, with the element farthest from $\chi$ being discarded from $\mathcal{N}$.
\end{itemize}

\smallskip

After performing $N_{\mathcal{S}}$ such iterations, we obtain the set $\mathcal{N}$ of $k$ nearest neighbors of $\chi$. Using the inverse distances from $\chi$ to the elements of $\mathcal{N}$ as weights (23), we assign a class to $\chi$.

\smallskip

We now estimate the computational cost of a single kNN iteration. The~feature vector $\chi(A)$, used to classify the matrix $A$, is 4-dimensional. Let $C_{\rho}$ denote the cost (in clock cycles) of computing the distance $\rho(\chi,\psi_i)$ in step 1. As the distance between vectors $x,y\!\in\!\mathbb{R}^4$, we use the squared Euclidean distance: $$\rho(x,y)=(x\!-\!y)^{\top}(x\!-\!y).$$ This avoids the computationally expensive square root operation required for~the standard Euclidean norm. On modern processors, a single-precision square root takes 10--16 clock cycles, while double-precision requires 12--20 cycles. By using the squared distance, we eliminate this cost entirely. 

The computation of $\rho(x,y)$ involves: 
\begin{itemize}
\item[--] 4 subtractions,
\item[--] 4 multiplications (for the element-wise products),
\item[--] 3 additions (to sum the products).
\end{itemize}

\noindent Assuming each addition and multiplication of single-precision numbers costs 4 clock cycles (a typical value for modern CPUs/GPUs), the total cost is: $$C_{\rho}=4C_++C_{\mathrm{sc}}(4)=4C_++4C_{\times}+3C_+=7C_++4C_{\times}=7\cdot 4+4\cdot 4=44.$$

The second step of the kNN iteration --- maintaining the set $\mathcal{N}$ of the $k$ nearest neighbors --- involves only comparisons and memory address updates. We neglect memory overhead and consider only comparison operations, which require 1--3 clock cycles per comparison (for both single- and double-precision values). In the worst case, up to $k$ comparisons are needed to insert a new candidate into $\mathcal{N}$. Thus, the total cost of one kNN iteration is bounded by: $$C_{\mathrm{iter}}=44+3\cdot k.\eqno(24)$$

From (18) and (24), it follows that to satisfy condition (14) when preconditioning the system, it suffices to enforce the inequality $$N_{\mathcal{S}}\cdot (44+3\cdot k)\le\alpha\cdot 4(14m+17n)\cdot k_{\mathrm{CG}^c}.$$ Consequently, the maximum admissible training sample satisfying (14) is bounded by $$N_{\mathcal{S}}\le\alpha\cdot\frac{4(14m+17n)\cdot k_{\mathrm{CG}^c}}{44+3\cdot k}=\alpha\cdot\frac{4(14c+17)n\cdot k_{\mathrm{CG}^c}}{44+3\cdot k}\eqno(25)$$ where $c\!=\! m/n$ denotes the average number of nonzeros per row.

In the absence of preconditioning, and using (19), the corresponding bound becomes $$N_{\mathcal{S}}\le\alpha\cdot\frac{4(6c+17)n\cdot k_{\mathrm{CG}}}{44+3\cdot k}\eqno(26).$$\\

\subsection{On the Inapplicability of Neural Networks\\ to the Prediction Problem} 

Neural networks might appear to be a promising approach for predicting the optimal tolerance $\hat{\varepsilon}_1$. Indeed, suppose that for an input matrix $A$, the value $\hat{\varepsilon}_1$ is obtained via a~function $\varphi(A)$, implemented by a neural network pretrained on a~dataset of feature vectors derived from a representative sample of matrices, i.e., $\hat{\varepsilon}_1\!=\!\varphi(A)$. In this scenario, the prediction cost would be $O(1)$, whereas using kNN incurs a cost of $O(N_{\mathcal{S}})$.

However, despite the low dimensionality of our feature vector, neural networks prove unsuitable for this prediction task. The reason lies in the high sensitivity of the target function: a small perturbation in the matrix $A$ --- that is, minor changes in the positions or values of its nonzero entries --- can lead to significant variations in the features $\tilde{\ell}(A)$ (pseudo-diameter) and $v(A)$ (residual decay rate) included in the vector $\chi(A)$. This sensitivity undermines the smoothness assumptions typically required for successful neural network training.

This conclusion is supported by extensive computational experiments using a multilayer perceptron (MLP) with the following hyperparameter ranges:

\smallskip

\begin{itemize}
\item[--] number of layers: from $3$ to $6$;
\item[--] neurons per layer: from $4$ to $30$;
\item[--] learning rates: from $0.0001$ to $0.1$ (step size $0.005$);
\item[--] batch sizes: from $1$ to $50$ (increment $10$);
\end{itemize}

\smallskip

\noindent (with training extended up to $10\, 000$ epochs). When trained on the same datasets used for kNN, the MLP achieved classification accuracy no better than that of kNN. More critically, when the data were split into training and test sets, the test accuracy dropped significantly, indicating poor generalization. This behavior confirms that the underlying mapping $A\to\hat{\varepsilon}_1$ lacks the regularity needed for effective neural network approximation.\\

\section{Computational Experiments}

\subsection{Samples used in computational experiments}

In the computational experiments, the results of which are presented and discussed below, samples of symmetric positive definite sparse matrices were used, constructed as follows.

Each matrix is derived from the adjacency matrix of a pre-generated unweighted undirected graph. Positive definiteness is ensured by setting the diagonal entries to $$a_{ii}=\mu\sum\limits_{j\neq i}|a_{ij}|,$$ where $\mu\!>\!1$ is a parameter controlling the condition number of the resulting matrix.

The sparsity pattern of each matrix is thus determined by the topology of its underlying graph. We consider the following graph types:

\smallskip

\begin{itemize}
\item[1.] {\bf Perturbed star graphs}: Random graphs generated by adding random edges to extended star graphs. Matrices derived from these graphs allow us to isolate and assess the influence of the matrix graph diameter on the efficiency of the proposed approach.

\smallskip

\item[2.] {\bf Connected random graphs:}  Graphs obtained by adding $cn$ random edges to random trees, where $c\!\in\!\{0.1, 0.5, 1, 2, 3, 5, 7, 10\}$  controls the edge density of the resulting graph.

\smallskip

\item[3.] {\bf Banded graphs:} Graphs corresponding to banded matrices of prescribed sparsity. Within a fixed band, each entry is nonzero independently with probability  $p\!\in\!\{0.4, 0.6, 0.8\}$. The band half-width is chosen to ensure sparsity even when the band is fully populated: for $n\!=\!500$, it does not exceed $25$; for $n\!=\!1000$, it does not exceed $50$.
 
\end{itemize}

\smallskip

When conducting experiments on a general sample $\mathcal{S}_0$ of feature vectors of size $N$, we form training $\mathcal{S}$ and test (validation) $\mathcal{T}$ samples, $\mathcal{S}\cap\mathcal{T}\!=\!\varnothing$. Elements of $\mathcal{S}$ and $\mathcal{T}$ are selected uniformly at random from the index set $\{1,\ldots,N\}$.

\smallskip

This procedure yields matrix feature vector samples for the following three matrix types:

\smallskip

\begin{itemize}
\item[1.] {\bf Extended star graph matrices with added random edges.} Matrix sizes: $n\!=\!301, 501, 1001$. Corresponding general sample sizes: $|\mathcal{S}_0|\!=\!3000$, $3600$, $4800$. 

\smallskip

For $n\!=\!1001$, an additional experiment was conducted by augmenting the original sample with $1600$ newly generated matrices. This demonstrates that increasing the representativeness of $\mathcal{S}_0$ --- in terms of the distribution of feature vectors used for matrix classification and $\varepsilon_1$ selection --- improves both classification accuracy and the efficiency of {\bf Algorithm~II} according to criterion (12).

\smallskip

\item[2.] {\bf Random sparse matrices.} Matrix sizes: $n\!=\!300, 500, 1000$. General sample size: $|\mathcal{S}_0|\!=\!4800$ for all $n$.

\smallskip

Each $n\!\times\! n$ matrix $A\!=\!(a_{ij})$ (as well as for banded matrices, see below) is constructed from the adjacency matrix of its underlying graph by replacing each nonzero entry (initially $1$) with a random value drawn from the interval $(-10, 10)$ using {\bf Algorithm~III} (see Appendix).

\smallskip

\item[3.] {\bf Random sparse banded matrices.} Matrix sizes: $n\!=\!300, 500, 1000$. General sample size: $|\mathcal{S}_0|\!=\!1680, 3000, 6000$, respectively.
\end{itemize}

\smallskip

The size of the general sample $\mathcal{S}_0$ for each matrix type is determined by the graph structure, the number of matrices generated for the given matrix and graph parameters, and the requirement that all resulting matrices remain sparse.

The training sample sizes for kNN classification --- computed via formula (26) with $\alpha\!=\!0.01$ and using the average number of CG iterations (over all matrices in $\mathcal{S}_0$) when executed entirely in double precision --- are relatively small. Specifically, they amount to at most $12\%$, $50\%$, and $9\%$ of the corresponding general sample sizes for the three matrix types, respectively. All elements of $\mathcal{S}_0$ not included in the training set form the test set.

\smallskip

Further details on the matrix sample generation procedure used in our computational experiments are provided in the Appendix.\\

\subsection{Efficiency criteria for the algorithm}

To compute average efficiency values of {\bf Algorithm~II} over the test sets --- and to compare them with the efficiency of {\bf Algorithm~I} executed using the optimal $\varepsilon_1$ for each input matrix --- we evaluated the following quantities for every matrix $A$ in the sample $\mathcal{S}_0^M$:

\smallskip

\begin{itemize}
\item[1)] $I_{\hat{\varepsilon}_1}(A)$: the value of $I(A)$ computed via (3) for {\bf Algorithm~I} with the optimal choice $\varepsilon_1\!=\!\hat{\varepsilon}_1(A)$;

\smallskip

\item[2)] $I(A)$: the value obtained by {\bf Algorithm~II} when selecting $\varepsilon_1$ via kNN using the matrix parameters $(n,m,\tilde{\ell},v)$;

\smallskip

\item[3)] $I_{\mathrm{dp}}(A)$: the number of iterations required by CG entirely in double precision to achieve the target tolerance $\varepsilon_2$.

\end{itemize}

\smallskip

For test samples $\mathcal{T}$ of size $N_{\mathcal{T}}$, we compute the total cost over all matrices in the sample: $$I_{\hat{\varepsilon}_1}(\mathcal{T})=\sum\limits_{i=1}^{N_{\mathcal{T}}}I_{\hat{\varepsilon}_1}(A_i),\qquad I(\mathcal{T})=\sum\limits_{i=1}^{N_{\mathcal{T}}}I(A_i),\qquad I_{\mathrm{dp}}(\mathcal{T})=\sum\limits_{i=1}^{N_{\mathcal{T}}}I_{\mathrm{dp}}(A_i).$$ 

The {\it overall efficiency} of {\bf Algo\-rithm I} (with the optimal $\varepsilon_1$) $E_{\mathrm{opt}}$ and {\bf Algo\-rithm I} --- relative to entirely double-precision CG --- is then defined as: $$E_{\mathrm{opt}}=\biggl (1-\frac{I_{\hat{\varepsilon}_1}(\mathcal{T})}{I_{\mathrm{dp}}(\mathcal{T})}\biggr )\cdot 100\%\mbox{\ \ and\ \ }E=\biggl (1-\frac{I(\mathcal{T})}{I_{\mathrm{dp}}(\mathcal{T})}\biggr )\cdot 100\%.\eqno(27)$$ 

We compare the efficiency values (27) for two choices of the parameter $\omega$:
\begin{itemize}
\item[(i)] the fixed estimate $\omega\!=\!1/3$ from (3)--(4), and
\item[(ii)] matrix-specific values of $\omega$ computed individually for each sample matrix by averaging over all iterations performed in single and double precision.
\end{itemize}

\subsection{Results of computational experiments}

This section summarizes the main findings of our computational experiments; complete results are provided in the Appendix.

\smallskip

\subsubsection{Classification accuracy}

For all matrix types considered, classification accuracy increases with matrix size, provided that the number of neighbors $k$ is appropriately chosen and the training set is sufficiently representative in terms of feature vector distribution (see the additional experiment on extended star graphs with randomly added edges). For matrices of size $n\!=\!1000$ (or $n\!=\!1001$ for graph-based matrices), the average classification accuracy is approximately \textbf{70\%}, \textbf{71\%}, and \textbf{74\%} for the three matrix types, respectively.

The standard deviation of the accuracy estimates decreases with increasing matrix size and, for $n\!=\!1000$ (or $1001$), does not exceed \textbf{3.3\%}, \textbf{1.3\%}, and \textbf{2.0\%}, respectively.

\subsubsection{Algorithm efficiency}

The efficiency of {\bf Algorithm II} also improves with matrix size. For $n\!=\!1000$ (or $1001$) and $\omega\!=\!1/3$, the average efficiency is approximately \textbf{22\%}, \textbf{22\%}, and \textbf{17.7\%} across the three matrix types. The corresponding standard deviations are at most \textbf{0.4\%}, \textbf{0.1\%}, and \textbf{0.4\%}.

Moreover, the performance gap between {\bf Algorithm~II} and {\bf Algorithm~I} (with optimal $\varepsilon_1$) is small: it does not exceed \textbf{3.0\%}, \textbf{1.5\%}, and \textbf{1.5\%}, respectively.
	
\smallskip 

We also performed computations using a matrix-specific value of $\omega$. In this setting, $\omega$ was determined from the ratio $t_{\mathrm{sp}}/t_{\mathrm{dp}}$, where $t_{\mathrm{sp}}$ and $t_{\mathrm{dp}}$ denote the average durations of a CG iteration for a given matrix in single and double precision, respectively. These timings were obtained individually for each matrix by averaging over all CG iterations executed during Algorithm I runs with all candidate values $\varepsilon_1\!=\!10^{-2},\ldots,10^{-7}$, as part of the matrix sample generation process for each matrix type and parameter configuration. All experiments were conducted on a Google Cloud server accessed via Google Colab using CPU.

For these matrix-specific $\omega$ values, the optimal $\varepsilon_1$ used in the efficiency evaluation of {\bf Algorithm~II} remained the same as in the $\omega\!=\!1/3$ case.

With this adaptive $\omega$, the efficiency of {\bf Algorithm~II} also increases with matrix size. For $n\!=\!1000$ (or $n\!=\!1001$ for graph-based matrices), the average efficiencies are \textbf{31.0\%}, \textbf{30.8\%}, and \textbf{25.3\%} for the three matrix types, respectively. The corresponding standard deviations do not exceed \textbf{0.50\%}, \textbf{0.15\%}, and \textbf{0.51\%}. The performance gap relative to {\bf Algorithm~I} (with optimal $\varepsilon_1$) remains small, at most \textbf{3.0\%}, \textbf{1.5\%}, and \textbf{1.5\%}, respectively.\\

\section*{Conclusions}

This paper presents an algorithm for optimal precision switching in the conjugate gradient (CG) method when solving sparse linear systems. Based on matrix parameters --- estimated in time which is negligible compared to the runtime of double-precision CG --- the algorithm predicts the optimal single-precision tolerance $\varepsilon_1$ that minimizes total mixed-precision computation time. For a given input matrix, this tolerance is determined by classifying its feature vector --- comprising matrix size $n$, number of nonzeros $m$, pseudo-diameter of the matrix graph, and the average CG convergence rate during early single-precision iterations --- using the $k$-nearest neighbors method.

We demonstrate that, in the absence of preconditioning, the diameter of the matrix graph governs computational efficiency for random sparse matrices with condition number $\kappa\!>\!10$ and random sparse banded matrices with $\kappa\!>\!100$.

The proposed algorithm reduces the effective computational cost --- measured in equivalent double-precision CG iterations --- by approximately \textbf{22\%} for random sparse matrices and \textbf{17.7\%} for banded matrices, assuming a $3\!:\!1$ double-to-single-precision iteration time ratio. When this ratio is estimated individually per matrix, the efficiency improves to \textbf{30.8\%} and \textbf{25.3\%}, respectively. In both cases, the performance loss relative to the oracle choice of $\varepsilon_1$ (i.e., with perfect knowledge of the optimal tolerance) does not exceed \textbf{1.5\%}.

Our experiments show that the algorithm's efficiency increases with matrix size.

%

%

\newpage

\section*{Appendix}

\section{Samples used in the experiments}

\subsection{Samples of extended star graph matrices\\ with added random edges}

A {\it path} is a graph consisting of a sequence of vertices $v_1,\ldots, v_{lc}$ connected by edges $(v_i,v_{i+1})$, $i\!=\!\overline{1,lc}$. An {\it extended star} is a graph with a single central vertex connected to $nc$ rays, each of which is a path of length $lc$, measured as the number of vertices in it.

We generated matrix samples from precomputed adjacency matrices of extended stars. The parameters of these stars are listed in Tables 7, 8, and 9, where:

\begin{itemize}
\item[--] $lc$ is the number of vertices in each ray,
\item[--] $nc$ is the number of rays.
\end{itemize}

For matrix sizes $n\!=\!301,501,1001$, the corresponding $(lc,nc)$ pairs satisfy $n\!=\!1+nc\cdot lc$ (accounting for the central vertex).

To each extended star, we added $\xi$ random edges, where $\xi$ is a random integer uniformly sampled from $[0,0.1\cdot n)$. The endpoints of each new edge $(i,j)$ were chosen uniformly at random from all unordered pairs $\{i,j\}$ with $1\!\le\!i\!<\!j\!\le\!n$, ensuring no self-loops or duplicate edges.

For each parameter pair $(lc,nc)$, we generated $100$ distinct matrices.

\begin{table}[H]
\begin{center}
\caption{Parameters $(lc,nc)$ of extended stars with added random edges for matrix sampling with $n\!=\!301$.}
\bigskip

\begin{tabular}{|c||c|c|c|c|c|c|c|c|c|c|}
\hline  
$lc$  & 1  & 2 & 5 & 10 & 15 & 20 & 30 & 60 &  150 & 300\\
\hline
$nc$  & 300 & 150 & 60 & 30 & 20 & 15 & 10 & 5 & 2 & 1\\
\hline  
\end{tabular}  
\end{center}  
\end{table}

\begin{table}[H]
\begin{center}
\caption{Parameters $(lc,nc)$ of extended stars with added random edges for matrix sampling with $n\!=\!501$.}

\bigskip

\begin{tabular}{|c||c|c|c|c|c|c|c|c|c|c|c|c|}
\hline  
$lc$  & 1  & 2 & 4 & 5 & 10 & 20 & 25 & 50 &  100 & 125 & 250 & 500\\
\hline
$nc$  & 500 & 250 & 125 & 100 & 50 & 25 & 20 & 10 & 5 & 4 & 2 & 1\\
\hline  
\end{tabular}  
\end{center}  
\end{table}

\begin{table}[H]
\begin{center}
\caption{Parameters $(lc,nc)$ of extended stars with added random edges for sampling matrices with $n\!=\!1001$.}

\bigskip

\begin{tabular}{|c||c|c|c|c|c|c|c|c|c|c|c|c|c|c|c|c|}
\hline  
$lc$  & 1  & 2 & 4 & 5 & 8 & 10 & 20 & 25 &  40 & 50 & 100 & 125 & 200 & 250 & 500 & 1000\\
\hline
$nc$  & 1000 & 500 & 250 & 200 & 125 & 100 & 50 & 40 & 25 & 20 & 10 & 8 & 5 & 4 & 2 & 1\\
\hline  
\end{tabular}  
\end{center}  
\end{table}

Using matrices with binary entries $a_{ij}\!\in\!\{0,1\}$ allows us to isolate the effect of the matrix graph structure on the performance of the proposed approach, independent of the influence of varying magnitudes in the matrix entries. The same matrices were used to illustrate:
\begin{itemize}
\item[(i)] the dependence of residual norms on the number of iterations in different precisions for star and path graph matrices (Figs.~2 and~3), and
\item[(ii)] the relationship between the efficiency of Algorithm I (with optimal $\varepsilon_1$) and the graph diameter (Table~1).
\end{itemize}

\smallskip

Note that, by choosing an appropriate parameter $\mu\!>\!1$, we can generate matrices of this type with arbitrarily large condition numbers. For all samples --- including the random and banded matrices discussed below --- the diagonal entries are defined as $$a_{ii}=\mu\cdot\sum\limits_{i\neq j} |a_{ij}|,$$ which ensures strict diagonal dominance. Three values of $\mu$ were used to construct the samples: $1.1$, $3$, and $10$. Varying $\mu$ yields a more uniform distribution of matrices across optimal $\varepsilon_1$ values, thereby increasing the difficulty of the classification task. Tables showing these distributions for each sample are provided below.\\

\subsection{Samples of random sparse matrices}

To generate a random sparse matrix of size $n$, we first constructed a graph consisting of a tree on $n$ vertices. For each density parameter $c$, we then added $\lfloor cn\rfloor$  random edges. The values $c\!\in\!\{0.1,0.5,1,2,3,5,7,10\}$ were used for each matrix size $n\!\in\!\{300,500,1000\}$.

\begin{codebox}
\Procname{$\proc{\bf Algorithm III}$}
\li \mbox{\bf Input:}\ $n\!\times\! n$ adjacency matrix $A_0=(a^0_{ij})$.\\ 

\li \For $i\leftarrow 1$ \To $n$ \Do
\li    \For $j\leftarrow i+1$ \To $n$ \Do
\li    \If $a^0_{ij}=1$:\Then
\li       \If $\xi_d>0.5$:\Then 
\li 						\If $\xi_s>0.5$:\Then 
\li									$a_{ij}\leftarrow +\xi_{m_1}$
\li					 		\Else
\li									$a_{ij}\leftarrow -\xi_{m_1}$
	  						\End
\li				\Else 
\li             \If $\xi_s>0.5$:\Then 
\li                	$a_{ij}\leftarrow +\xi_{m_2}$
\li   					\Else
\li 						  	$a_{ij}\leftarrow -\xi_{m_2}$
    						\End
          \End 					
\li       $a_{ji}\leftarrow a_{ij}$.\\																					
        \End
        \End				
    \End
		
\li \mbox{\bf Output:}\ $n\!\times\! n$ matrix $A\!=\!(a_{ij})$.		
		
\end{codebox}

To generate samples of random and random banded matrices, we first constructed the adjacency matrix $A_0$ of the underlying sparsity graph and then produced the numerical matrix using {\bf Algorithm III}.

The matrix entries are generated as follows:
\begin{itemize}
\item[---] the random variable $\xi_d\!\sim\!\mathcal{U}(0,1)$ determines the magnitude range of the entry;
\item[---] the random variable $\xi_s\!\sim\!\mathcal{U}(0,1)$ determines its sign;
\item[---] the random variables $\xi_{m_1}\!\sim\!\mathcal{U} (a_1,b_1)$, $\xi_{m_2}\!\sim\!\mathcal{U}(a_2,b_2)$ determine the absolute value of the entry.
\end{itemize}

\noindent For all matrices in the random and banded samples considered below, we set $a_1\!=\!0$, $b_1\!=\!3$ and $a_2\!=\!7$, $b_2\!=\!10$. The intervals $[0,3]$ and $[7,10]$ for the absolute values of the matrix entries were chosen to ensure a broad and representative coverage of satisfactorily conditioned and moderately ill-conditioned matrices, according to our condition-number-based classification.\\

\subsection{Sparse banded matrix sampling}

We generated banded matrix graphs with full bandwidths $b\!\in\!\{3,\ldots,29\}$, $b\!\in\!\{3,\ldots,51\}$, and $b\!\in\!\{3,\ldots,101\}$ for matrix sizes $n\!=\!300$, $n\!=\!500$, and $n\!=\!1000$, respectively. The number of nonzero entries within the band was controlled by a density parameter $p\!\in\![0.4,0.6,0.8]$: for each $p$, edges inside the band were included independently with probability $p$. The numerical matrix entries were then generated from the adjacency matrix $A_0$ of the resulting graph using {\bf Algorithm~III}. For each pair $(p,b)$, we produced $20$ distinct matrices.

In cases where the graph was disconnected, the pseudo-diameter was defined as the maximum pseudo-diameter over all connected components.

Notably, banded matrix graphs typically yield well-conditioned matrices under our classification --- even for $p\!>\!0.5$. This favorable conditioning is attributed to their structural properties: banded graphs contain relatively few terminal vertices (i.e., vertices of degree one). In contrast, matrices whose graphs have many terminal vertices --- such as those derived from small-diameter extended stars --- tend to be ill-conditioned.\\

\section{Formation of training and test samples}

After constructing the general matrix samples of the types described above, we partitioned each into training and test samples. Let $N$ denote the size of the general sample of matrices $\mathcal{S}_0$. The size $N_{\mathcal{S}}$ of the training sample $\mathcal{S}$ is computed based on the parameters of formula (26): $$N_{\mathcal{S}}=\alpha\cdot\frac{4(6c+17)n\cdot k_{\mathrm{CG}}}{44+3\cdot k},$$ where
\begin{itemize} 
\item[1)] $\alpha\!=\!0.01$;
\item[2)] $k_{\mathrm{CG}}$  is the average number of CG iterations (over all matrices of the given type and size in $\mathcal{S}_0$ required to achieve a residual tolerance of $10^{-10}$ when executed entirely in double precision;
\item[3)] $c$  is set to the smallest value used in matrix generation, i.e., $c\!=\!0.1$;
\item[4)] the kNN parameter takes five values: $k\!\in\!\{1,5,10,15,20\}$.
\end{itemize}

Once $N_{\mathcal{S}}$ is determined, the test sample size is defined as $$N_{\mathcal{T}}=N-N_{\mathcal{S}},$$ so that $\mathcal{T}\!=\!\mathcal{S}_0\setminus \mathcal{S}$; that is, all matrices of a given type and size not selected for training are assigned to the test set. Consequently, the training set is substantially smaller than the test set. The exact sizes of $\mathcal{S}$ and $\mathcal{T}$ are reported alongside the corresponding experimental results.\\

\section{Choice of the parameter $\omega$\\ based on single-to-double precision speedup}

The experimental results for the efficiency of {\bf Algorithm~II} are presented below for two choices of the parameter $\omega$ in (3)--(4).

In the first variant, we fix $\omega\!=\!1/3$. In the second, $\omega$ is computed individually for each matrix as the empirical ratio $\omega=t_{\mathrm{sp}}/t_{\mathrm{dp}}$, where $t_{\mathrm{sp}}$ and $t_{\mathrm{dp}}$ denote the average durations of a CG iteration in single and double precision, respectively. These timings were obtained by averaging over all CG iterations performed for a given matrix across all experiments with $\varepsilon_1\!\in\!\{\!10^{-2},\ldots,10^{-7}\}$, as part of the matrix sample generation process for each matrix type and parameter configuration. All computations were carried out on a Google Cloud server accessed via Google Colab using CPU.

For both variants, the optimal $\varepsilon_1$ used in the efficiency evaluation of {\bf Algorithm~II} remains the same.\\

\section{Experiments}

The data in all tables below were obtained by averaging over 100 experiments with each resulting general sample. For each experiment, the general sample was split into training and test sets.

For every matrix sample considered, the following quantities are provided:
\begin{itemize}
\item[1)] the size $N\!=\!|\mathcal{S}_0|$ of the general sample, along with the sizes $N_{\mathcal{S}}\!=\!|\mathcal{S}|$ and $N_{\mathcal{T}}\!=\!|\mathcal{T}|$ of the training and test samples;
\item[2)] The distribution of matrices by optimal tolerance $\hat{\varepsilon}_1$: each entry in row $P$ and column $\hat{\varepsilon}_1$ gives the proportion $P$ (as a percentage) of matrices in $\mathcal{S}_0$ for which the optimal tolerance equals the column's $\hat{\varepsilon}_1$;
\item[3)] $k_{\mathrm{CG}}$: the average number of CG iterations (over all matrices in $\mathcal{S}_0$) when executed entirely in double precision;
\item[4)] classification accuracy $\mathsf{Acc}$ of the kNN classifier used in {\bf Algorithm~II};
\item[5)] localization measure $\mathsf{Loc}$ of feature vectors (defined as the average fraction of same-class neighbors among the $k$ nearest neighbors);
\item[6)] efficiency $E$ of {\bf Algorithm II} relative to full double-precision CG, reported for two choices of $\omega$: fixed ($\omega\!=\!1/3$) and adaptive ($\omega\!=\!t_{\mathrm{sp}}/t_{\mathrm{dp}}$).
\end{itemize}

All reported values are expressed as percentages. Columns prefixed with {\sffamily M} (e.g., {\sffamily M}Acc) denote mean values, while those prefixed with $\sigma$ (e.g., $\sigma\mbox{Acc}$) report the corresponding standard deviations. Efficiency values $E$ and $E_{\mathrm{opt}}$ from (29) are presented as $\mbox{\sffamily M}E/\mbox{\sffamily M}E_{\mathrm{opt}}$.

The proportions $P$ in the $\hat{\varepsilon}_1$ distribution tables are rounded to one decimal place ($0.1\%$), while all other results are rounded to two decimal places ($0.01\%$).

Finally, we note that a $3$-dimensional feature vector $(m, \tilde{\ell}, v)$ --- excluding matrix size $n$ --- was used in the experiments to demonstrate that the efficiency of {\bf Algorithm~II} increases with $n$.\\

\subsection{Extended star graph matrices with random edges added}

Sample parameters for extended star graph matrices. 

\subsubsection{Sizes of training and test samples}

For $n\!=\!301$: $N\!=\!3000$, $k_{\mathrm{CG}}\!=\!30.34$.

\begin{table}[H]
\begin{center}
\caption{Sizes of training and test samples for different values of the kNN parameter $k$. $n\!=\!301$.}

\bigskip

\begin{tabular}{|c|c|c|c|c|c|}
\hline  
$k$ & $1$ & $5$  & $10$  & $15$ & $20$\\  
\hline
$N_{\mathcal{S}}$ & 137 & 109 & 87 & 65 & 62\\  
\hline  
$N_{\mathcal{T}}$ & 2863 & 2891 & 2913 & 2935 & 2938\\  
\hline  
\end{tabular}  
\end{center}  
\end{table}

For $n\!=\!501$: $N\!=\!3600$, $k_{\mathrm{CG}}\!=\!31.5$.

\begin{table}[H]
\begin{center}
\caption{Sizes of training and test samples for different values of the kNN parameter $k$. $n\!=\!501$.}

\bigskip

\begin{tabular}{|c|c|c|c|c|c|}
\hline  
$k$ & $1$ & $5$  & $10$  & $15$ & $20$\\  
\hline
$N_{\mathcal{S}}$ & 236 & 188 & 150 & 112 & 107\\  
\hline  
$N_{\mathcal{T}}$ & 3364 & 3412 & 3450 & 3488 & 3493\\  
\hline  
\end{tabular}  
\end{center}  
\end{table}

For $n\!=\!1001$: $N\!=\!4800$, $k_{\mathrm{CG}}\!=\!34.46$.

\begin{table}[H]
\begin{center}
\caption{Sizes of training and test samples for different values of the kNN parameter $k$. $n\!=\!1001$.}

\bigskip

\begin{tabular}{|c|c|c|c|c|c|}
\hline  
$k$ & $1$ & $5$  & $10$  & $15$ & $20$\\  
\hline
$N_{\mathcal{S}}$ & 517 & 412 & 329 & 246 & 234\\  
\hline  
$N_{\mathcal{T}}$ & 4283 & 4388 & 4471 & 4554 & 4566\\  
\hline  
\end{tabular}  
\end{center}  
\end{table}

\subsubsection{Distributions of matrices from samples by the value $\hat{\varepsilon}_1$}

\begin{table}[H]
\begin{center}
\caption{Distribution of matrices by the value $\hat{\varepsilon}_1$. $n\!=\!301$.}

\bigskip

\begin{tabular}{|c|c|c|c|c|c|c|}
\hline  
$\hat{\varepsilon}_1$ & $10^{-2}$ & $10^{-3}$  & $10^{-4}$  & $10^{-5}$ & $10^{-6}$ & $10^{-7}$\\  
\hline  
$P$ & 16.7 & 8.0 & 25.4 & 43.2 & 6.6 & 0.0\\  
\hline  
\end{tabular}  
\end{center}  
\end{table}

\begin{table}[H]
\begin{center}
\caption{Distribution of matrices by the value $\hat{\varepsilon}_1$. $n\!=\!501$.}

\bigskip

\begin{tabular}{|c|c|c|c|c|c|c|}
\hline  
$\hat{\varepsilon}_1$ & $10^{-2}$ & $10^{-3}$  & $10^{-4}$  & $10^{-5}$ & $10^{-6}$ & $10^{-7}$\\  
\hline  
$P$ & 17.9 & 14.0 & 33.3 & 34.3 & 0.5 & 0.0\\
\hline  
\end{tabular}  
\end{center}  
\end{table}

\begin{table}[H]
\begin{center}
\caption{Distribution of matrices by the value $\hat{\varepsilon}_1$. $n\!=\!1001$.}

\bigskip

\begin{tabular}{|c|c|c|c|c|c|c|}
\hline  
$\hat{\varepsilon}_1$ & $10^{-2}$ & $10^{-3}$  & $10^{-4}$  & $10^{-5}$ & $10^{-6}$ & $10^{-7}$\\  
\hline  
$P$ & 18.8 & 16.6 & 36.0 & 28.4 & 0.1 & 0.0 \\  
\hline  
\end{tabular}  
\end{center}  
\end{table}

\subsubsection{kNN accuracy and localization of feature vectors}

\begin{table}[H]
\begin{center}
\caption{Accuracy of kNN.}

\bigskip

\begin{tabular}{|c||c|c||c|c||c|c|}
\hline  
& \multicolumn{2}{|c||}{301}&\multicolumn{2}{|c||}{501}&\multicolumn{2}{|c|}{1001}\\  
\hline  
$k$ & $\mbox{\sffamily M}\mbox{Acc}$ & $\sigma\mbox{Acc}$ & $\mbox{\sffamily M}\mbox{Acc}$ & $\sigma\mbox{Acc}$ & $\mbox{\sffamily M}\mbox{Acc}$ & $\sigma\mbox{Acc}$\\  
\hline  
1  & 67.92 & 3.24 & 68.74 & 2.41 & 65.23 & 1.78\\  
\hline  
5  & 68.68 & 4.49 & 68.71 & 3.19 & 66.47 & 2.17\\  
\hline  
10 & 66.36 & 4.40 & 67.96 & 3.30 & 67.04 & 2.56\\  
\hline  
15 & 63.85 & 7.18 & 68.30 & 5.16 & 66.24 & 3.10\\  
\hline  	
20 & 62.89 & 8.80 & 66.03 & 6.29 & 65.56 & 3.64\\  
\hline  
\end{tabular}  
\end{center}  
\end{table}

\begin{table}[H]
\begin{center}
\caption{Localization of feature vectors.}

\bigskip

\begin{tabular}{|c||c|c||c|c||c|c|}
\hline  
& \multicolumn{2}{|c||}{301}&\multicolumn{2}{|c||}{501}&\multicolumn{2}{|c|}{1001}\\  
\hline  
$k$ & $\mbox{\sffamily M}\mbox{Loc}$ & $\sigma\mbox{Loc}$ & $\mbox{\sffamily M}\mbox{Loc}$ & $\sigma\mbox{Loc}$ & $\mbox{\sffamily M}\mbox{Loc}$ & $\sigma\mbox{Loc}$\\  
\hline  
1  & 67.92 & 3.24 & 68.74 & 2.41 & 65.23 & 1.78\\  
\hline  
5  & 67.17 & 3.34 & 67.11 & 2.35 & 65.17 & 1.81\\    
\hline  
10 & 66.50 & 3.34 & 66.15 & 2.07 & 64.95 & 1.82\\  
\hline  
15 & 67.41 & 4.30 & 66.66 & 3.04 & 64.62 & 2.13\\  
\hline  	
20 & 67.64 & 4.54 & 65.81 & 3.62 & 64.53 & 2.27\\  
\hline  
\end{tabular}  
\end{center}  
\end{table}

\paragraph{Improving kNN accuracy without increasing the training sample size.} To the general sample of feature vectors from extended star matrices with added random edges ($n\!=\!1001$), we added $1600$ feature vectors from newly generated matrices of the same type, using $\mu\!=\!1.1$ and the graph parameters specified in Table~9. This augmentation significantly improves kNN classification accuracy, while the localization measure changes only marginally (see Tables 19 and 20). Consequently, the efficiency of {\bf Algorithm~II} increases for both choices of $\omega$ (see Tables 23 and 24). The training set size $N_{\mathcal{T}}$ remains unchanged (as listed in Table 12), while the test set size $N_{\mathcal{T}}$ increases by $1600$.

\begin{table}[H]
\begin{center}

\caption{Distribution of matrices by the value of $\hat{\varepsilon}_1$, after augmenting the general sample with 1600 additional matrices (cf. Table~15).}

\bigskip

\begin{tabular}{|c|c|c|c|c|c|c|}
\hline  
$\hat{\varepsilon}_1$ & $10^{-2}$ & $10^{-3}$  & $10^{-4}$  & $10^{-5}$ & $10^{-6}$ & $10^{-7}$\\  
\hline  
$P$ & 15.2 & 12.8 & 27.8 & 43.5 & 0.7 & 0.0 \\
\hline  
\end{tabular}  
\end{center}  
\end{table}

\begin{table}[H]
\begin{center}
\caption{Accuracy of kNN, after augmenting the general sample with 1600 additional matrices (cf. Table 16).}

\bigskip

\begin{tabular}{|c||c|c||c|c||c|c|}
\hline  
& \multicolumn{2}{|c||}{301}&\multicolumn{2}{|c||}{501}&\multicolumn{2}{|c|}{1001}\\  
\hline  
$k$ & $\mbox{\sffamily M}\mbox{Acc}$ & $\sigma\mbox{Acc}$ & $\mbox{\sffamily M}\mbox{Acc}$ & $\sigma\mbox{Acc}$ & $\mbox{\sffamily M}\mbox{Acc}$ & $\sigma\mbox{Acc}$\\  
\hline  
1  & 67.92 & 3.24 & 68.74 & 2.41 & 67.79 & 1.72\\  
\hline  
5  & 68.68 & 4.49 & 68.71 & 3.19 & 70.15 & 1.96\\  
\hline  
10 & 66.36 & 4.40 & 67.96 & 3.30 & 72.07 & 2.50\\  
\hline  
15 & 63.85 & 7.18 & 68.30 & 5.16 & 73.83 & 2.85\\  
\hline  	
20 & 62.89 & 8.80 & 66.03 & 6.29 & 75.04 & 3.22\\  
\hline  
\end{tabular}  
\end{center}  
\end{table}

\begin{table}[H]
\begin{center}
\caption{Localization of feature vectors, after augmenting the general sample with $1600$ additional matrices (cf.  Table 17).}

\bigskip

\begin{tabular}{|c||c|c||c|c||c|c|}
\hline  
& \multicolumn{2}{|c||}{301}&\multicolumn{2}{|c||}{501}&\multicolumn{2}{|c|}{1001}\\  
\hline  
$k$ & $\mbox{\sffamily M}\mbox{Loc}$ & $\sigma\mbox{Loc}$ & $\mbox{\sffamily M}\mbox{Loc}$ & $\sigma\mbox{Loc}$ & $\mbox{\sffamily M}\mbox{Loc}$ & $\sigma\mbox{Loc}$\\  
\hline  
1  & 67.92 & 3.24 & 68.74 & 2.41 & 67.79 & 1.72\\  
\hline  
5  & 67.17 & 3.34 & 67.11 & 2.35 & 66.88 & 1.53\\    
\hline  
10 & 66.50 & 3.34 & 66.15 & 2.07 & 66.51 & 1.85\\  
\hline  
15 & 67.41 & 4.30 & 66.66 & 3.04 & 66.32 & 2.21\\  
\hline  	
20 & 67.64 & 4.54 & 65.81 & 3.62 & 66.56 & 2.31\\  
\hline  
\end{tabular}  
\end{center}  
\end{table}

\subsubsection{Efficiency of the algorithm}

\begin{table}[H]
\begin{center}
\caption{Efficiency $E$ of the algorithm relative to double-precision CG for~$w\!=\!1/3$.}

\bigskip

\begin{tabular}{|c||c|c||c|c||c|c|}
\hline  
& \multicolumn{2}{|c||}{301}&\multicolumn{2}{|c||}{501}&\multicolumn{2}{|c|}{1001}\\  
\hline  
$k$ & $\mbox{\sffamily M}E$ & $\sigma E$ & $\mbox{\sffamily M}E$ & $\sigma E$ & $\mbox{\sffamily M}E$ & $\sigma E$\\  
\hline  
1  & 21.88/24.62 & 0.35/0.05 & 19.10/21.87 & 0.28/0.06 & 18.65/21.86 & 0.20/0.07\\  
\hline  
5  & 21.87/24.62 & 0.04/0.48 & 18.86/21.87 & 0.35/0.06 & 18.71/22.86 & 0.24/0.06\\  
\hline  
10 & 21.66/24.62 & 0.57/0.04 & 18.79/21.88 & 0.39/0.04 & 18.76/21.85 & 0.30/0.06\\  
\hline  
15 & 21.54/24.61 & 0.92/0.03 & 18.81/21.87 & 0.56/0.05 & 18.66/21.75 & 0.41/0.05\\  
\hline  	
20 & 21.52/24.62 & 1.04/0.04 & 18.58/21.88 & 0.72/0.04 & 18.58/21.85 & 0.48/0.05\\  
\hline  
\end{tabular}  
\end{center}  
\end{table}

\begin{table}[H]
\begin{center}
\caption{Efficiency $E$ of the algorithm relative to double-precision CG for~$w\!=\!t_{\mathrm{sp}}/t_{\mathrm{dp}}$.}

\bigskip

\begin{tabular}{|c||c|c||c|c||c|c|}
\hline  
& \multicolumn{2}{|c||}{301}&\multicolumn{2}{|c||}{501}&\multicolumn{2}{|c|}{1001}\\  
\hline  
$k$ & $\mbox{\sffamily M}E$ & $\sigma E$ & $\mbox{\sffamily M}E$ & $\sigma E$ & $\mbox{\sffamily M}E$ & $\sigma E$\\  
\hline  
1  & 14.57/16.92 & 0.25/0.05 & 24.22/27.15 & 0.33/0.05 & 27.50/31.81 & 0.27/0.06\\  
\hline  
5  & 14.46/16.91 & 0.33/0.05 & 23.93/27.15 & 0.45/0.05 & 27.55/31.21 & 0.36/0.05\\  
\hline  
10 & 14.49/16.92 & 0.45/0.04 & 23.93/29.13 & 0.52/0.05 & 27.51/31.22 & 0.37/0.05\\  
\hline  
15 & 14.34/16.89 & 0.67/0.03 & 23.75/29.13 & 0.66/0.05 & 27.37/31.21 & 0.53/0.03\\  
\hline  	
20 & 14.16/16.90 & 0.88/0.03 & 23.69/27.14 & 0.87/0.05 & 27.41/27.20 & 0.68/0.04\\  
\hline  
\end{tabular}  
\end{center}  
\end{table}

\paragraph{Improving the accuracy of kNN without increasing the training set size.} Tables 23 and 24 present results for matrix samples of the same type as above. The general sample for $n\!=\!1001$ is augmented with additional matrices to improve kNN accuracy.

\begin{table}[H]
\begin{center}
\caption{Efficiency $E$ of the algorithm relative to double-precision CG for~$\omega\!=\!1/3$, after augmenting the general sample with $1600$ additional matrices (cf. Table~21).}

\bigskip

\begin{tabular}{|c||c|c||c|c||c|c|}
\hline  
& \multicolumn{2}{|c||}{300}&\multicolumn{2}{|c||}{500}&\multicolumn{2}{|c|}{1000}\\  
\hline  
$k$ & $\mbox{\sffamily M}E$ & $\sigma E$ & $\mbox{\sffamily M}E$ & $\sigma E$ & $\mbox{\sffamily M}E$ & $\sigma E$\\  
\hline  
1  & 16.20/17.57 & 0.14/0.05 & 21.74/22.99 & 0.09/0.04 & 21.76/24.4 & 0.18/0.02\\  
\hline  
5  & 16.28/17.56 & 0.12/0.05 & 21.87/23.00 & 0.09/0.03 & 21.88/24.4 & 0.24/0.03\\  
\hline  
10 & 16.30/17.57 & 0.14/0.05 & 21.90/22.99 & 0.11/0.03 & 22.07/24.4 & 0.26/0.02\\  
\hline  
15 & 16.29/17.56 & 0.18/0.05 & 21.90/23.00 & 0.13/0.02 & 22.20/24.40 & 0.30/0.02\\  
\hline  	
20 & 16.32/17.56 & 0.24/0.05 & 21.92/22.99 & 0.13/0.02 & 22.26/24.40 & 0.31/0.02\\  
\hline  
\end{tabular}  
\end{center}  
\end{table}

\begin{table}[H]
\begin{center}
\caption{Efficiency $E$ of the algorithm relative to double-precision CG for~$\omega\!=\!t_{\mathrm{sp}}/t_{\mathrm{dp}}$, after augmenting the general sample with $1600$ additional matrices (cf. Table~22).}

\bigskip

\begin{tabular}{|c||c|c||c|c||c|c|}
\hline  
& \multicolumn{2}{|c||}{300}&\multicolumn{2}{|c||}{500}&\multicolumn{2}{|c|}{1000}\\  
\hline  
$k$ & $\mbox{\sffamily M}E$ & $\sigma E$ & $\mbox{\sffamily M}E$ & $\sigma E$ & $\mbox{\sffamily M}E$ & $\sigma E$\\  
\hline  
1  & 7.56/8.83 & 0.09/0.05 & 24.10/25.36 & 0.10/0.05 & 30.72/33.90 & 0.25/0.03\\  
\hline  
5  & 7.60/8.82 & 0.09/0.04 & 24.22/25.35 & 0.11/0.05 & 30.95/33.92 & 0.28/0.04\\  
\hline  
10 & 7.62/8.81 & 0.11/0.03 & 24.26/25.35 & 0.12/0.05 & 31.18/33.92 & 0.38/0.04\\  
\hline  
15 & 7.62/8.81 & 0.12/0.03 & 24.27/25.35 & 0.16/0.05 & 31.41/33.96 & 0.40/0.05\\  
\hline  	
20 & 7.64/8.80 & 0.14/0.02 & 24.28/25.36 & 0.15/0.05 & 31.55/33.97 & 0.43/0.05\\  
\hline  
\end{tabular}  
\end{center}  
\end{table}

\subsection{Random sparse matrices}

Sample parameters for random matrices.

\subsubsection{Sizes of training and test samples}

For $n\!=\!300$: $N\!=\!4800$, $k_{\mathrm{CG}}\!=\!71.92$.

\begin{table}[H]
\begin{center}
\caption{Sizes of training and test samples for different values of the kNN parameter $k$. $n\!=\!300$.}

\bigskip

\begin{tabular}{|c|c|c|c|c|c|}
\hline  
$k$ & $1$ & $5$  & $10$  & $15$ & $20$\\  
\hline
$N_{\mathcal{S}}$ & 323 & 258 & 205 & 153 & 146\\  
\hline  
$N_{\mathcal{T}}$ & 4477 & 4542 & 4595 & 4647 & 4654\\  
\hline  
\end{tabular}  
\end{center}  
\end{table}

For $n\!=\!500$: $N\!=\!4800$, $k_{\mathrm{CG}}\!=\!87.57$.

\begin{table}[H]
\begin{center}
\caption{Sizes of training and test samples for different values of the kNN parameter $k$. $n\!=\!500$.}

\bigskip

\begin{tabular}{|c|c|c|c|c|c|}
\hline  
$k$ & $1$ & $5$  & $10$  & $15$ & $20$\\  
\hline
$N_{\mathcal{S}}$ & 656 & 532 & 417 & 311 & 296\\  
\hline  
$N_{\mathcal{T}}$ & 4144 & 4268  & 4383 & 4489 & 4504 \\  
\hline  
\end{tabular}  
\end{center}  
\end{table}

For $n\!=\!1000$: $N\!=\!4800$, $k_{\mathrm{CG}}\!=\!105.17$.

\begin{table}[H]
\begin{center}
\caption{Sizes of training and test samples for different values of the kNN parameter $k$. $n\!=\!1000$.}

\bigskip

\begin{tabular}{|c|c|c|c|c|c|}
\hline  
$k$ & $1$ & $5$  & $10$  & $15$ & $20$\\  
\hline
$N_{\mathcal{S}}$ & 1575 & 1255 & 1001 & 749 & 712\\    
\hline  
$N_{\mathcal{T}}$ & 3225 & 3545 & 3799 & 4051 & 4088\\   
\hline  
\end{tabular}  
\end{center}  
\end{table}

\subsubsection{Distributions of matrices from samples by the value $\hat{\varepsilon}_1$}

\begin{table}[H]
\begin{center}
\caption{Distribution of matrices by the value $\hat{\varepsilon}_1$. $n\!=\!300$.}

\bigskip

\begin{tabular}{|c|c|c|c|c|c|c|}
\hline  
$\hat{\varepsilon}_1$ & $10^{-2}$ & $10^{-3}$  & $10^{-4}$  & $10^{-5}$ & $10^{-6}$ & $10^{-7}$\\  
\hline  
$P$& 5.6 & 23.5 & 43.8 & 26.8 & 0.3 & 0.0\\
\hline  
\end{tabular}  
\end{center}  
\end{table}

\begin{table}[H]
\begin{center}
\caption{Distribution of matrices by the value $\hat{\varepsilon}_1$. $n\!=\!500$.}

\bigskip

\begin{tabular}{|c|c|c|c|c|c|c|}
\hline  
$\hat{\varepsilon}_1$ & $10^{-2}$ & $10^{-3}$  & $10^{-4}$  & $10^{-5}$ & $10^{-6}$ & $10^{-7}$\\  
\hline  
$P$ & 5.9 & 24.5 & 41.1 & 28.1 & 0.4 & 0.0\\
\hline  
\end{tabular}  
\end{center}  
\end{table}

\begin{table}[H]
\begin{center}
\caption{Distribution of matrices by the value $\hat{\varepsilon}_1$. $n\!=\!1000$.}

\bigskip

\begin{tabular}{|c|c|c|c|c|c|c|}
\hline  
$\hat{\varepsilon}_1$ & $10^{-2}$ & $10^{-3}$  & $10^{-4}$  & $10^{-5}$ & $10^{-6}$ & $10^{-7}$\\  
\hline  
$P$& 5.6 & 25.8 & 36.8 & 31.5 & 0.3 & 0.0\\
\hline  
\end{tabular}  
\end{center}  
\end{table}

\subsubsection{kNN accuracy and localization of feature vectors}

\begin{table}[H]
\begin{center}
\caption{Accuracy of kNN.}

\bigskip

\begin{tabular}{|c||c|c||c|c||c|c|}
\hline  
& \multicolumn{2}{|c||}{300}&\multicolumn{2}{|c||}{500}&\multicolumn{2}{|c|}{1000}\\  
\hline  
$k$ & $\mbox{\sffamily M}\mbox{Acc}$ & $\sigma\mbox{Acc}$ & $\mbox{\sffamily M}\mbox{Acc}$ & $\sigma\mbox{Acc}$ & $\mbox{\sffamily M}\mbox{Acc}$ & $\sigma\mbox{Acc}$\\  
\hline  
1  & 69.64 & 2.17 & 69.23 & 1.54 & 68.93 & 1.22\\  
\hline  
5  & 70.44 & 2.18 & 70.86 & 1.70 & 70.35 & 1.29\\  
\hline  
10 & 70.40 & 2.59 & 70.91 & 2.10 & 71.07 & 1.39\\  
\hline  
15 & 69.79 & 3.81 & 70.41 & 2.59 & 71.07 & 1.71\\  
\hline  	
20 & 70.34 & 4.60 & 70.58 & 2.71 & 71.57 & 1.91\\  
\hline  
\end{tabular}  
\end{center}  
\end{table}

\begin{table}[H]
\begin{center}
\caption{Localization of feature vectors.}

\bigskip

\begin{tabular}{|c||c|c||c|c||c|c|}
\hline  
& \multicolumn{2}{|c||}{300}&\multicolumn{2}{|c||}{500}&\multicolumn{2}{|c|}{1000}\\  
\hline  
$k$ & $\mbox{\sffamily M}\mbox{Loc}$ & $\sigma\mbox{Loc}$ & $\mbox{\sffamily M}\mbox{Loc}$ & $\sigma\mbox{Loc}$ & $\mbox{\sffamily M}\mbox{Loc}$ & $\sigma\mbox{Loc}$\\  
\hline  
1  & 69.64 & 2.17 & 69.23 & 1.54 & 68.93 & 1.22\\  
\hline  
5  & 69.78 & 1.80 & 69.66 & 1.30 & 68.90 & 0.96\\  
\hline  
10 & 70.40 & 2.59 & 69.46 & 1.56 & 69.07 & 1.15\\  
\hline  
15 & 70.61 & 2.39 & 69.48 & 1.73 & 69.17 & 1.13\\  
\hline  	
20 & 71.04 & 2.57 & 70.00 & 1.70 & 69.45 & 1.25\\  
\hline  	
\end{tabular}  
\end{center}  
\end{table}

\subsubsection{Efficiency of the algorithm}

\begin{table}[H]
\begin{center}
\caption{Efficiency $E$ of the algorithm relative to double-precision CG for~$\omega\!=\!1/3$.}

\bigskip

\begin{tabular}{|c||c|c||c|c||c|c|}
\hline  
& \multicolumn{2}{|c||}{300}&\multicolumn{2}{|c||}{500}&\multicolumn{2}{|c|}{1000}\\  
\hline  
$k$ & $\mbox{\sffamily M}E$ & $\sigma E$ & $\mbox{\sffamily M}E$ & $\sigma E$ & $\mbox{\sffamily M}E$ & $\sigma E$\\  
\hline  
1  & 16.20/17.57 & 0.14/0.05 & 21.74/22.99 & 0.09/0.04 & 21.91/23.13 & 0.08/0.06\\  
\hline  
5  & 16.28/17.56 & 0.12/0.05 & 21.87/23.00 & 0.09/0.03 & 22.04/23.14 & 0.08/0.05\\  
\hline  
10 & 16.30/17.57 & 0.14/0.05 & 21.90/22.99 & 0.11/0.03 & 22.11/23.14 & 0.07/0.06\\  
\hline  
15 & 16.29/17.56 & 0.18/0.05 & 21.90/23.00 & 0.13/0.02 & 22.12/23.14 & 0.10/0.05\\  
\hline  	
20 & 16.32/17.56 & 0.24/0.05 & 21.92/22.99 & 0.13/0.02 & 22.14/23.12 & 0.10/0.04\\  
\hline  
\end{tabular}  
\end{center}  
\end{table}

\begin{table}[H]
\begin{center}
\caption{Efficiency $E$ of the algorithm relative to double-precision CG for~$\omega\!=\!t_{\mathrm{sp}}/t_{\mathrm{dp}}$.}

\bigskip

\begin{tabular}{|c||c|c||c|c||c|c|}
\hline  
& \multicolumn{2}{|c||}{300}&\multicolumn{2}{|c||}{500}&\multicolumn{2}{|c|}{1000}\\  
\hline  
$k$ & $\mbox{\sffamily M}E$ & $\sigma E$ & $\mbox{\sffamily M}E$ & $\sigma E$ & $\mbox{\sffamily M}E$ & $\sigma E$\\  
\hline  
1  & 7.56/8.83 & 0.09/0.05 & 24.10/25.36 & 0.10/0.05 & 30.56/31.86 & 0.09/0.06\\  
\hline  
5  & 7.60/8.82 & 0.09/0.04 & 24.22/25.35 & 0.11/0.05 & 30.82/31.87 & 0.12/0.05\\  
\hline  
10 & 7.62/8.81 & 0.11/0.03 & 24.26/25.35 & 0.12/0.05 & 30.83/31.85 & 0.12/0.05\\  
\hline  
15 & 7.62/8.81 & 0.12/0.03 & 24.27/25.35 & 0.16/0.05 & 30.83/31.80 & 0.15/0.04\\  
\hline  	
20 & 7.64/8.80 & 0.14/0.02 & 24.28/25.36 & 0.15/0.05 & 30.82/31.79 & 0.13/0.04\\  
\hline  
\end{tabular}  
\end{center}  
\end{table}

\subsection{Sparse banded matrices}

Sampling parameters for sparse banded matrices.

\subsubsection{Sizes of training and test samples}

For $n\!=\!301$: $N\!=\!1680$, $k_{\mathrm{CG}}\!=\!41.41$.

\begin{table}[H]
\begin{center}
\caption{Sizes of training and test samples for different values of the kNN parameter $k$. $n\!=\!300$.}

\bigskip

\begin{tabular}{|c|c|c|c|c|c|}
\hline  
$k$ & $1$ & $5$  & $10$  & $15$ & $20$\\  
\hline
$N_{\mathcal{S}}$ & 186 & 148 & 118 & 88 & 84\\  
\hline  
$N_{\mathcal{T}}$ & 1494 & 1532 & 1562 & 1592 & 1596\\  
\hline  
\end{tabular}  
\end{center}  
\end{table}

For $n\!=\!500$: $N\!=\!3000$, $k_{\mathrm{CG}}\!=\!36.66$.

\begin{table}[H]
\begin{center}
\caption{Sizes of training and test samples for different values of the kNN parameter $k$. $n\!=\!500$.}

\bigskip

\begin{tabular}{|c|c|c|c|c|c|}
\hline  
$k$ & $1$ & $5$  & $10$  & $15$ & $20$\\  
\hline
$N_{\mathcal{S}}$ & 275 & 219 & 174 & 130 & 124\\  
\hline  
$N_{\mathcal{T}}$ & 2725 & 2781 & 2826 & 2870 & 2876\\  
\hline  
\end{tabular}  
\end{center}  
\end{table}

For $n\!=\!1000$: $N\!=\!6000$, $k_{\mathrm{CG}}\!=\!33.81$.

\begin{table}[H]
\begin{center}
\caption{Sizes of training and test samples for different values of the kNN parameter $k$. $n\!=\!1000$.}

\bigskip

\begin{tabular}{|c|c|c|c|c|c|}
\hline  
$k$ & $1$ & $5$  & $10$  & $15$ & $20$\\  
\hline
$N_{\mathcal{S}}$ & 506 & 403 & 322 & 240 & 229\\    
\hline  
$N_{\mathcal{T}}$ & 5494 & 5597 & 5678 & 5760 & 5771\\   
\hline  
\end{tabular}  
\end{center}  
\end{table}

\subsubsection{Distributions of matrices from samples by the value $\hat{\varepsilon}_1$}

\begin{table}[H]
\begin{center}
\caption{Distribution of matrices by the value $\hat{\varepsilon}_1$. $n\!=\!300$.}

\bigskip

\begin{tabular}{|c|c|c|c|c|c|c|}
\hline  
$\hat{\varepsilon}_1$ & $10^{-2}$ & $10^{-3}$  & $10^{-4}$  & $10^{-5}$ & $10^{-6}$ & $10^{-7}$\\  
\hline  
$P$ & 24.9 & 40.4 & 26.8 & 7.5 & 0.4 & 0.0\\  
\hline  
\end{tabular}  
\end{center}  
\end{table}

\begin{table}[H]
\begin{center}
\caption{Distribution of matrices by the value $\hat{\varepsilon}_1$. $n\!=\!500$.}

\bigskip

\begin{tabular}{|c|c|c|c|c|c|c|}
\hline  
$\hat{\varepsilon}_1$ & $10^{-2}$ & $10^{-3}$  & $10^{-4}$  & $10^{-5}$ & $10^{-6}$ & $10^{-7}$\\  
\hline  
$P$ & 41.9 & 36.8 & 16.4 & 4.5 & 0.4 & 0.0\\
\hline  
\end{tabular}  
\end{center}  
\end{table}

\begin{table}[H]
\begin{center}
\caption{Distribution of matrices by the value $\hat{\varepsilon}_1$. $n\!=\!1000$.}

\bigskip

\begin{tabular}{|c|c|c|c|c|c|c|}
\hline  
$\hat{\varepsilon}_1$ & $10^{-2}$ & $10^{-3}$  & $10^{-4}$  & $10^{-5}$ & $10^{-6}$ & $10^{-7}$\\  
\hline  
$P$ & 51.8 & 37.4 & 8 & 2.5 & 0.3 & 0.0\\
\hline  
\end{tabular}  
\end{center}  
\end{table}

\subsubsection{kNN accuracy and localization of feature vectors}

\begin{table}[H]
\begin{center}
\caption{Accuracy of kNN.}

\bigskip

\begin{tabular}{|c||c|c||c|c||c|c|}
\hline  
& \multicolumn{2}{|c||}{300}&\multicolumn{2}{|c||}{500}&\multicolumn{2}{|c|}{1000}\\  
\hline  
$k$ & $\mbox{\sffamily M}\mbox{Acc}$ & $\sigma\mbox{Acc}$ & $\mbox{\sffamily M}\mbox{Acc}$ & $\sigma\mbox{Acc}$ & $\mbox{\sffamily M}\mbox{Acc}$ & $\sigma\mbox{Acc}$\\  
\hline  
1  & 69.67 & 2.66 & 71.33 & 1.86 & 74.33 & 1.37\\  
\hline  
5  & 68.48 & 3.39 & 71.00 & 2.40 & 73.52 & 1.55\\  
\hline  
10 & 68.79 & 4.02 & 71.11 & 2.43 & 73.16 & 1.85\\  
\hline  
15 & 67.27 & 4.47 & 70.07 & 3.18 & 73.44 & 1.99\\  
\hline  	
20 & 67.50 & 5.60 & 70.49 & 3.17 & 73.08 & 1.86\\  
\hline  
\end{tabular}  
\end{center}  
\end{table}

\begin{table}[H]
\begin{center}
\caption{Localization of feature vectors.}

\bigskip

\begin{tabular}{|c||c|c||c|c||c|c|}
\hline  
& \multicolumn{2}{|c||}{300}&\multicolumn{2}{|c||}{500}&\multicolumn{2}{|c|}{1000}\\  
\hline  
$k$ & $\mbox{\sffamily M}\mbox{Loc}$ & $\sigma\mbox{Loc}$ & $\mbox{\sffamily M}\mbox{Loc}$ & $\sigma\mbox{Loc}$ & $\mbox{\sffamily M}\mbox{Loc}$ & $\sigma\mbox{Loc}$\\  
\hline  
1  & 69.67 & 2.66 & 71.33 & 1.86 & 74.33 & 1.37\\  
\hline  
5  & 68.51 & 2.60 & 71.21 & 1.73 & 74.39 & 1.15\\  
\hline  
10 & 68.79 & 4.02 & 70.75 & 1.88 & 73.74 & 1.35\\  
\hline  
15 & 67.22 & 3.43 & 69.58 & 2.42 & 73.14 & 1.39\\  
\hline  	
20 & 67.20 & 2.98 & 69.44 & 2.21 & 72.92 & 1.33\\  
\hline  	
\end{tabular}  
\end{center}  
\end{table}

\subsubsection{Efficiency of the algorithm}

\begin{table}[H]
\begin{center}
\caption{Efficiency $E$ of the algorithm relative to double-precision CG for~$\omega\!=\!1/3$.}

\bigskip

\begin{tabular}{|c||c|c||c|c||c|c|}
\hline  
& \multicolumn{2}{|c||}{300}&\multicolumn{2}{|c||}{500}&\multicolumn{2}{|c|}{1000}\\  
\hline  
$k$ & $\mbox{\sffamily M}E$ & $\sigma E$ & $\mbox{\sffamily M}E$ & $\sigma E$ & $\mbox{\sffamily M}E$ & $\sigma E$\\  
\hline  
1  & 15.48/16.97 & 0.26/0.09 & 15.74/17.13 & 0.20/0.09 & 17.77/19.02 & 0.24/0.06\\  
\hline  
5  & 15.47/16.98 & 0.31/0.08 & 15.72/17.08 & 0.22/0.10 & 17.78/19.01 & 0.26/0.06\\  
\hline  
10 & 15.49/16.95 & 0.32/0.08 & 15.59/17.01 & 0.21/0.09 & 17.70/19.02 & 0.29/0.05\\  
\hline  
15 & 15.37/16.96 & 0.36/0.06 & 15.67/17.11 & 0.23/0.11 & 17.70/19.02 & 0.30/0.04\\  
\hline  	
20 & 15.45/16.95 & 0.41/0.06 & 15.66/17.17 & 0.28/0.11 & 17.67/19.02 & 0.34/0.04\\  
\hline  
\end{tabular}  
\end{center}  
\end{table}

\begin{table}[H]
\begin{center}
\caption{Efficiency $E$ of the algorithm relative to double-precision CG for~$\omega\!=\!t_{\mathrm{sp}}/t_{\mathrm{dp}}$.}

\bigskip

\begin{tabular}{|c||c|c||c|c||c|c|}
\hline  
& \multicolumn{2}{|c||}{300}&\multicolumn{2}{|c||}{500}&\multicolumn{2}{|c|}{1000}\\  
\hline  
$k$ & $\mbox{\sffamily M}E$ & $\sigma E$ & $\mbox{\sffamily M}E$ & $\sigma E$ & $\mbox{\sffamily M}E$ & $\sigma E$\\  
\hline  
1  & 8.61/10.06 & 0.20/0.12 & 15.74/17.13 & 0.20/0.09 & 25.52/26.76 & 0.35/0.06\\  
\hline  
5  & 8.49/9.92 & 0.22/0.09  & 15.72/17.08 & 0.22/0.10 & 25.37/26.61 & 0.37/0.05\\  
\hline  
10 & 8.38/9.91 & 0.22/0.07  & 15.59/17.01 & 0.21/0.09 & 25.35/26.69 & 0.46/0.05\\  
\hline  
15 & 8.48/10.00 & 0.26/0.11 & 15.67/17.11 & 0.23/0.11 & 25.35/26.69 & 0.49/0.05\\  
\hline  	
20 & 8.45/10.02 & 0.28/0.14 & 15.66/17.17 & 0.28/0.11 & 25.26/26.71 & 0.51/0.04\\  
\hline  
\end{tabular}  
\end{center}  
\end{table}


\begin{thebibliography}{5}
%
\bibitem{Higham}
Higham, N.J., Mary, T.: \textit{Mixed precision algorithms in numerical linear algebra} // Acta Numerica. Vol.~31, 347--414 (2022).

\bibitem{Abdelfattah} Abdelfattah, A., Anzt, H., Boman, E.G., Carson, E., Cojean, T., Dongarra, J.,  Fox, A., Gates, M., Higham, N.J., Li, X.S., Loe, J., Luszczek, P., Pranesh, S., Rajamanickam, S., Ribizel, T., Smith, B.F., Swirydowicz, K., Thomas, S., Tomov, S., Tsai, Y.M. and Yang, U.M.: \textit{A survey of numerical linear algebra methods utilizing mixed-precision arithmetic} // Int. J. High Perform. Comput. Appl. Vol.~35, 344--369 (2021).

\bibitem{Turner} 
Turner, K., Walker, H.F.: \textit{Efficient high accuracy solutions with GMRES(m)} // SIAM J. Sci. Statist. Comput. Vol. 12, 815--825 (1992).

\bibitem{Lindqwist} 
Lindquist, N., Luszczek, P., Dongarra, J.: \textit{Improving the performance of the GMRES method using mixed-precision techniques}// Communications in Computer and Information Science (J. Nichols et al., eds), Springer, pp. 51--66 (2020).

\bibitem{Giraud}
Giraud, L., Haidar, A., Watson, L.T.: \textit{Mixed-precision preconditioners in parallel domain decomposition solvers} // Domain Decomposition Methods in Science and Engineering XVII (U. Langer et al., eds). Vol. 60 of Lecture Notes in Computational Science and Engineering, Springer, 357--364 (2008).

\bibitem{Gratton}
Gratton, S., Simon, E., Titley-Peloquin, D., Toint, P.: \textit{Exploiting variable precision in GMRES}. Available at arXiv:1907.10550. (2019).

\bibitem{Aliaga}
Aliaga, J.I., Anzt, H., Grutzmacher, T., Quintana-Orti, E.S., Tomas, A.E.: \textit{Compressed basis GMRES on high performance GPUs}. Available at arXiv:2009.12101. (2020).


\bibitem{Saad2003}
Saad, Y. Iterative Methods for Sparse Linear Systems. 2nd~ed. Philadelphia: SIAM, 2003.

\bibitem{Greenbaum1997}
Greenbaum, A. Iterative Methods for Solving Linear Systems. Philadelphia: SIAM, 1997.

\bibitem{Simoncini2007}
Simoncini, V., Szyld, D.~B. \textit{Recent computational developments in Krylov subspace methods for linear systems} // Numer. Linear Algebra Appl. Vol.~14, no.~1, 1--59 (2007).

\bibitem{Chow2015}
Chow,~E., Patel,~A. \textit{Fine-grained parallel incomplete LU factorization} // SIAM J. Sci. Comput. Vol.~37, no.~2. C169--C193 (2015).

\bibitem{Rude2026} 
R\"ude, U. \textit{Conjugate gradient methods are not efficient: experimental study of the locality limitation}. https://arxiv.org/pdf/2601.10322 (2026)

\bibitem{Huawei} Parameter optimization for restarted mixed precision iterative sparse solver. Conference MOTOR-2024 Challenges. https://motor24.oscsbras.ru/challenges/Challenge3.pdf

\bibitem{Carson2024} 
Carson, E., Liesen, J., Strako\v{s}, Z. \textit{Towards understanding CG and GMRES through examples} // Linear Algebra and It's Applications. Vol.~692, no.~164, 241--291 (2024).

\bibitem{VanDamHaemers1995} E.R. Van Dam, W.H. Haemers. \textit{Eigenvalues and the Diameter of Graphs} // Linear and Multilinear Algebra. Vol.~39, 33--44, (1995).
DOI: 10.1080/03081089508818378


\bibitem{Burkhardt} 
Burkhardt, P. \textit{Optimal algebraic breadth-first search for sparse graphs} // ACM Transactions on Knowledge Discovery from Data (TKDD). Vol.~15(5). Art. No.: 77, 1--19. (2021) https://doi.org/10.1145/3446216

\end{thebibliography}
\end{document}